\documentclass[11pt]{article}
\usepackage[a4paper,margin=1in]{geometry}
\usepackage[T1]{fontenc}
\usepackage{amsmath,amssymb,amsthm,mathtools}
\usepackage{booktabs}
\usepackage{array}
\usepackage{enumitem}
\usepackage{microtype}
\usepackage{lmodern}
\usepackage{url}
\usepackage{xurl}
\usepackage{aliascnt}
\usepackage[colorlinks=true,linkcolor=blue,citecolor=blue,urlcolor=blue]{hyperref}
\usepackage[nameinlink,capitalise,noabbrev]{cleveref}

\newtheorem{theorem}{Theorem}[section]
\newaliascnt{lemma}{theorem}
\newtheorem{lemma}[lemma]{Lemma}
\aliascntresetthe{lemma}
\newaliascnt{proposition}{theorem}
\newtheorem{proposition}[proposition]{Proposition}
\aliascntresetthe{proposition}
\newaliascnt{corollary}{theorem}
\newtheorem{corollary}[corollary]{Corollary}
\aliascntresetthe{corollary}
\newaliascnt{definition}{theorem}

\aliascntresetthe{definition}

\theoremstyle{remark}
\newaliascnt{remark}{theorem}
\newtheorem{remark}[remark]{Remark}
\aliascntresetthe{remark}
\crefname{lemma}{Lemma}{Lemmas}
\Crefname{lemma}{Lemma}{Lemmas}
\crefname{corollary}{Corollary}{Corollaries}
\Crefname{corollary}{Corollary}{Corollaries}
\crefname{claim}{Claim}{Claims}
\Crefname{claim}{Claim}{Claims}

\title{\textbf{A Fan-Type Theorem for Hamiltonian Berge Cycles in 3-Uniform Hypergraphs}}
\author{Ailian Chen\thanks{Corresponding author.}\\
School of Mathematics and Statistics, Fuzhou University, Fuzhou, Fujian Province 350108, China\\
\texttt{elian1425@fzu.edu.cn}
\and
Yanqiang Liu\\
School of Mathematics and Statistics, Fuzhou University, Fuzhou, Fujian Province 350108, China\\
\texttt{lyq.0217@qq.com}}
\date{}

\begin{document}
\maketitle

\begin{abstract}
Let $H$ be a simple $3$-uniform hypergraph on $n\ge7$ vertices, and put
$\tau_n=\binom{\lfloor(n-1)/2\rfloor}{2}$. We prove that $H$ has a
Hamiltonian Berge cycle if its $2$-shadow is $2$-connected, every pair
at Berge distance two has an endpoint of hypergraph degree at least
$\tau_n+1$, and every pair of codegree one has both endpoints above
$\tau_n$. The distance-two threshold is sharp, and the order range
cannot be extended to $n=4,5,6$. The proof separates two structural
arguments. A low-degree endpoint of an extremal longest Berge path
forces Hamiltonicity, independently of the order and the degree
threshold. For high-degree endpoints, injective charging and its
equality cases force pairs of codegree one that permit further
rotations. The cases $n=7,8$ are settled by exhaustive verification,
with a completeness proof and a second implementation using direct
vertex--hyperedge cycle search.
\end{abstract}

\medskip
{\small

}

\section{Introduction}\label{sec:introduction}

Fan's theorem states that a $2$-connected graph $G$ on $n$ vertices is
Hamiltonian whenever
\[
 d_G(u,v)=2 \quad\Longrightarrow\quad
 \max\{d_G(u),d_G(v)\}\ge n/2.
\]
Unlike the degree conditions of Dirac~\cite{1}, Ore~\cite{2}, and
P\'osa~\cite{3}, Fan's condition concerns only pairs at distance two~\cite{4}.

A Berge cycle consists of distinct vertices and distinct hyperedges
representing consecutive vertex pairs~\cite{5,6}. Kostochka, Luo, and
McCourt~\cite{7} determined the sharp minimum-degree threshold for
Hamiltonian Berge cycles in uniform hypergraphs. For $3$-graphs this
threshold is $\binom{\lfloor(n-1)/2\rfloor}{2}+1$. P\'osa- and Ore-type
conditions have also been studied~\cite{8,9,13}; see~\cite{10,11,12}
for related results. Recent work of Wang, Gerbner, and Zhao~\cite{15}
studies extremal and stability questions under a minimum-degree hypothesis.
That stability setting is distinct from the present local distance-two
condition, which admits vertices far below the Dirac threshold.

We ask whether the same threshold can be imposed locally, on an endpoint
of each pair at distance two. Our theorem allows vertices of degree
three at every order $n\ge7$. The additional hypothesis requires
codegree at least two on shadow edges incident with vertices below the
threshold. Pairs joining two high-degree vertices may still have a
unique representative. This is precisely where an argument in the
shadow alone ceases to suffice: a Hamilton cycle of the shadow need
not admit distinct representatives. We control these unique-representative
pairs by charging hyperedges at the endpoints of a longest Berge path.

All hypergraphs are finite, simple, and $3$-uniform; we write
``$3$-graph'' for short. The $2$-shadow $\partial_2H$ has vertex set
$V(H)$, with two vertices adjacent whenever they lie in a common
hyperedge.  For distinct $u,v$ write
\[
 E_H(u,v)=\{e\in E(H):\{u,v\}\subseteq e\},\qquad
 \deg_H(u,v)=|E_H(u,v)|.
\]
A Berge path or cycle is an alternating sequence of pairwise distinct key
vertices and pairwise distinct representative hyperedges, each representative
containing its two adjacent key vertices.  The Berge distance $d_H(u,v)$ is
the minimum length of a Berge path joining $u$ and $v$.  A hypergraph is connected precisely when its $2$-shadow is connected.
We call $H$ $2$-edge-connected if $H-e$ is connected for every
$e\in E(H)$. When $H$ is fixed, we omit the subscript $H$ from
$\deg_H(u)$, $\deg_H(u,v)$, $E_H(u,v)$, and $N_H(u)$.  We retain $d_H$
for Berge distance to distinguish it from ordinary graph distance.

By Lemma~\ref{G:lem:shadow-distance}, Berge distance two agrees with
graph distance two in the shadow. The additional condition concerns
shadow edges with a unique representative.

\begin{theorem}[Main theorem]\label{global:main}
Let $H$ be a simple $3$-uniform hypergraph on $n\ge7$ vertices
such that $\partial_2H$ is $2$-connected.  Put
\[
 \tau_n=\binom{\lfloor(n-1)/2\rfloor}{2}.
\]
Assume that for every pair of distinct vertices $u,v$,
\begin{enumerate}[label=\textup{(\alph*)}]
\item if the Berge distance $d_H(u,v)=2$, then
\[
 \max\{\deg_H(u),\deg_H(v)\}\ge \tau_n+1;
\]
\item if the codegree $\deg_H(u,v)=1$, then
\[
 \min\{\deg_H(u),\deg_H(v)\}\ge \tau_n+1.
\]
\end{enumerate}
Then $H$ contains a Hamiltonian Berge cycle.
\end{theorem}

\begin{remark}[Equivalent multiplicity form]\label{global:mult-form}
Condition~\ref{global:main}\textup{(b)} is equivalent to the following
statement:
\begin{equation}\label{global:low-mult}
 \deg_H(x)\le \tau_n\ \text{ and }\ xy\in E(\partial_2H)
 \quad\Longrightarrow\quad
 \deg_H(x,y)\ge2.
\end{equation}
Indeed, a positive codegree at a sub-threshold vertex cannot equal one;
conversely, \eqref{global:low-mult} excludes a low endpoint of a codegree-one pair.
\end{remark}

\subsection*{Sharpness and examples}
The theorem applies beyond the minimum-degree setting. Let $W$ be a set
of $n-1$ vertices, choose distinct $a,b,c\in W$, take all triples in
$W$, and add a vertex $x$ with the three incident triples
$\{x,a,b\}$, $\{x,a,c\}$, and $\{x,b,c\}$. Then $\deg_H(x)=3$,
while every vertex in $W$ has degree greater than $\tau_n$. The shadow is
$K_{n-1}$ together with the three edges joining $x$ to $a,b,c$.
Each shadow edge incident with $x$ has codegree two, so
Theorem~\ref{global:main} applies although $\deg_H(x)<\tau_n+1$.

\begin{remark}[Sharpness of the distance-two threshold]\label{L:rem:sharp-a}
We use the standard extremal construction of Kostochka, Luo, and
McCourt~\cite{7} (Construction~2 in the cited preprint), specialized to
$3$-graphs. Let $m=\lfloor(n-1)/2\rfloor$, let $A$ and $B$ be disjoint with
$|A|=m$ and $|B|=n-m$, and let $H_1(n)$ consist of all triples contained in
$A$ together with all triples having two vertices in $A$ and one in $B$.
For every $b\in B$,
\[
 \deg_{H_1}(b)=\binom m2=\tau_n,
\]
while every two vertices of $B$ are at Berge distance two.  No pair has
codegree one, and
\[
 \partial_2H_1=K_m\vee\overline K_{n-m}
\]
is $2$-connected.  The hypergraph is also $2$-edge-connected. Since every
hyperedge meets $B$ in at most one vertex, consecutive key vertices of
a Berge cycle cannot both lie in $B$.  Since $|B|>|A|$, no
Hamiltonian Berge cycle exists.  Hence the threshold $\tau_n+1$ in
Theorem~\ref{global:main}\textup{(a)} cannot be lowered.
\end{remark}

\begin{remark}[Why codegree-one control is necessary]\label{L:rem:necess-b}
Let $H_2$ have vertex set $\{u_1,\ldots,u_n\}$, contain all triples on
$\{u_3,\ldots,u_n\}$, and in addition contain precisely
$\{u_1,u_2,u_3\}$ and $\{u_1,u_2,u_4\}$.  Then $H_2$ is
$2$-edge-connected, its $2$-shadow is $2$-connected, and the distance-two
condition of Theorem~\ref{global:main}\textup{(a)} holds.  But $u_1$ and
$u_2$ are both incident only with the same two hyperedges.  If they are
nonconsecutive in a Hamiltonian Berge cycle, the four cycle incidences at
$u_1,u_2$ require four distinct representatives chosen from those two edges;
if they are consecutive, the three incidences meeting them still require
three distinct representatives.  Neither is possible, so condition~\textup{(b)} cannot be omitted.
By Remark~\ref{global:mult-form}, it requires at least two representatives
for each shadow edge incident with a sub-threshold vertex. Whether this
requirement can be weakened remains open.
\end{remark}

\begin{proposition}[Sharpness of the order range]\label{global:smallcounter}
For each $n\in\{4,5,6\}$ there is a simple $2$-edge-connected
$3$-uniform hypergraph $H_n$ whose $2$-shadow is $2$-connected, which
satisfies the two local hypotheses of Theorem~\ref{global:main} with the
corresponding threshold, but has no Hamiltonian Berge cycle.
\end{proposition}
\begin{proof}
Let
\[
 V(H_n)=\{x,v_1,\ldots,v_{n-1}\},\qquad
 E(H_n)=\{\{x,v_i,v_{i+1}\}:i\in\mathbb Z_{n-1}\},
\]
with rim indices read cyclically.  This cone over $C_{n-1}$ has $n-1$ hyperedges, but a Hamiltonian Berge
cycle requires $n$ distinct representatives.  Its shadow is a wheel and is $2$-connected, and deleting
any hyperedge leaves a connected hypergraph.  The centre has degree $n-1$
and every rim vertex degree two.  For $n=5,6$ one has $\tau_n+1=2$, and for
$n=4$ one has $\tau_n+1=1$.  Every rim pair $v_iv_{i+1}$ has codegree one
and every pair $xv_i$ has codegree two, so both local hypotheses hold.
\end{proof}

Put $G=\partial_2H$ and $r=\lfloor(n-1)/2\rfloor$. Since $H$ is simple,
\begin{equation}\label{eq:shadow-vs-hyper-degree}
\deg_H(v)\le \binom{\deg_G(v)}2\qquad(v\in V(H)).
\end{equation}
Consequently
\[
 \deg_H(v)\ge\binom r2+1\quad\Longrightarrow\quad
 \deg_G(v)\ge r+1\ge n/2.
\]
By Lemma~\ref{G:lem:shadow-distance}, Theorem~\ref{global:main}(a)
ensures that every two vertices at distance two in $G$ have a vertex
of shadow degree at least $n/2$.  Fan's theorem gives a Hamilton cycle in $G$, but this cycle need not lift
to a Berge cycle: distinct cycle edges may have no distinct representative
hyperedges available. The following lemma gives a sufficient condition
for such a lift.

\begin{lemma}[A safe lifting criterion]\label{intro:hall-lifting}
Let $H$ be a simple $3$-graph on $n\ge4$ vertices, and let
$C=v_1v_2\cdots v_nv_1$ be a Hamilton cycle of $\partial_2H$.
If $\deg_H(v_i,v_{i+1})\ge2$ for every $i$ (indices modulo $n$),
then $H$ has a Hamiltonian Berge cycle with key-vertex order $C$.
\end{lemma}
\begin{proof}
For each $i$ let $\mathcal F_i=E_H(v_i,v_{i+1})$.  Every
$\mathcal F_i$ contains at least two distinct hyperedges.  Since $n\ge4$,
a triple contains at most two edges of the fixed simple cycle $C$.
Thus for every index set $I\subseteq[n]$,
\[
2|I|\le\sum_{i\in I}|\mathcal F_i|
\le 2\left|\bigcup_{i\in I}\mathcal F_i\right|.
\]
Hall's theorem yields distinct representatives $e_i\in\mathcal F_i$.
Thus
\[
v_1,e_1,v_2,e_2,\ldots,v_n,e_n,v_1
\]
is a Hamiltonian Berge cycle.
\end{proof}

Condition~\ref{global:main}(b) guarantees two representatives only
when a shadow edge is incident with a vertex of degree at most
$\tau_n$. A shadow Hamilton cycle may therefore contain codegree-one
edges between high-degree vertices.

The codegree condition also makes a separate $2$-edge-connectivity
assumption unnecessary.
\begin{lemma}[Redundancy of $2$-edge-connectivity]\label{intro:edge-redundancy}
Let $H$ be a connected simple $3$-graph with $2$-connected $2$-shadow.
If $\delta(H)\ge2$, then $H-e$ is connected for every $e\in E(H)$.
Consequently the hypotheses of Theorem~\ref{global:main} imply that
$H$ is $2$-edge-connected.
\end{lemma}
\begin{proof}
Suppose first that $H-e$ is disconnected, with components
$X_1,\ldots,X_s$, $s\ge2$.  Every component meets $e$, since $H$ is
connected and the only deleted hyperedge is $e$.  As $|e|=3$, some
component $X$ intersects $e$ in precisely one vertex $x$.
If $X=\{x\}$, then $x$ belongs to no hyperedge of $H-e$, and
$\deg_H(x)=1$, contrary to the minimum degree.  If $|X|>1$, then
$G-x$ has no edge between $X\setminus\{x\}$ and $V(H)\setminus X$,
where $G=\partial_2H$: a hyperedge other than $e$ cannot meet distinct
components of $H-e$, and $e\cap X=\{x\}$.  These two sets are nonempty,
so $x$ is a cut vertex of $G$, a contradiction.

For the final assertion, $2$-connectivity of the shadow ensures
that every vertex has positive hypergraph degree.  If a vertex
$x$ had degree one, its unique edge $\{x,y,z\}$ would satisfy
$\deg_H(x,y)=1$ while $\deg_H(x)=1<\tau_n+1$, contradicting
Theorem~\ref{global:main}(b).  Thus $\delta(H)\ge2$, and the first
part applies.
\end{proof}

For comparison, Salia's P\'osa-type condition~\cite{8} for $3$-graphs of order $n>6$
requires, with $d_1\le\cdots\le d_n$ denoting the vertex degrees,
\[
 d_i>i\ \ (i=1,2),\qquad
 d_i>\binom i2\ \ \Bigl(3\le i\le\bigl\lfloor\tfrac{n-1}2\bigr\rfloor\Bigr),
\]
with the additional even-order requirement
$ d_{(n-2)/2}>\binom{(n-2)/2}2+1$.
The following example satisfies our hypotheses but not the latter
degree-sequence condition.

\begin{proposition}[Examples outside the P\'osa bound]\label{intro:posa-separation}
For $n=8$ and for every $n\ge10$, there is a simple $3$-graph satisfying
Theorem~\ref{global:main}(a)--(b) with $2$-connected shadow whose degree
sequence fails Salia's $3$-uniform P\'osa-type criterion~\cite{8}.
\end{proposition}
\begin{proof}
For $n=8$, take $V=\{x,y,a,b,c,d,e,f\}$ and include all triples of
$\{x,y,a,b\}$, all triples of $\{a,b,d,e,f\}$, and
\[
 \{c,a,b\},\quad\{c,a,d\},\quad\{c,a,e\},\quad\{c,b,d\}.
\]
The degrees of $x,y,c$ are $3,3,4$, respectively, and all other degrees
are at least six. The only low vertices $x,y$ are adjacent in the shadow,
and every shadow pair incident with either has codegree two. The shadow
is $2$-connected: the cliques on $\{x,y,a,b\}$ and $\{a,b,d,e,f\}$
meet in $ab$, and $c$ has neighbours $a,b,d,e$. Thus (a)--(b) hold,
whereas $d_3=4$ violates the even-order requirement $d_3>4$.

For $n\ge10$, let $L=\{1,2,3,4\}$ and let $W$ be a disjoint set of
$n-4$ vertices containing $a,b$. On $L\cup\{a,b\}$ take the ten triples
\[
 \mathcal D=\{123,124,13a,14b,1ab,23b,24a,2ab,34a,34b\},
\]
where juxtaposition denotes a set of three vertices, and put
$E(H)=\mathcal D\cup\binom W3$. Every pair of $L\cup\{a,b\}$
lies in exactly two triples of $\mathcal D$, and each of these six
vertices lies in five. Hence every vertex of $L$ has degree five;
vertices of $W\setminus\{a,b\}$ have degree $\binom{n-5}{2}$;
and $a,b$ have degree $5+\binom{n-5}{2}$. Writing
$r=\lfloor(n-1)/2\rfloor$, we have $n-5\ge r+1$ for $n\ge10$,
so all vertices of $W$ have degree greater than $\tau_n$, while
the four vertices of $L$ are low.

The shadow is the union of the cliques on $L\cup\{a,b\}$ and $W$,
which meet in $ab$, so it is $2$-connected. The low vertices form
a clique, giving (a); every positive codegree is at least two,
giving (b). The first four degrees are $5,5,5,5$, so
$d_4=5\le\binom42=6$. Since $4\le\lfloor(n-1)/2\rfloor$,
Salia's condition fails at index four for every $n\ge10$.
\end{proof}

These examples also fail the Dirac-type minimum-degree hypothesis.
The Ore-type theorem in~\cite{9} instead concerns degree sums in
the $2$-shadows of $[3]$-graphs, whose hyperedges have size at most
three. The examples establish separation from the stated P\'osa
criterion; no general non-comparability with all Ore-type conditions
is asserted.

\subsection*{Position of the result and proof outline}
The minimum-degree hypothesis $\delta(H)\ge\tau_n+1$ implies
(a)--(b). By \eqref{eq:shadow-vs-hyper-degree}, it also implies
$\delta(\partial_2H)\ge n/2$ and therefore $2$-connectivity of the
shadow. Our theorem thus extends the $3$-uniform Hamiltonian case
of~\cite{7} at the same sharp threshold. Proposition~\ref{intro:posa-separation}
shows that its scope also extends beyond the stated P\'osa condition
at arbitrarily large orders.

Condition (b) supplies two representatives only on shadow edges
incident with a low vertex. Requiring this on every shadow edge would
make Hall's theorem sufficient once Fan's theorem has produced a
shadow Hamilton cycle, by Lemma~\ref{intro:hall-lifting}. In the present
setting, representatives of high--high pairs can be unique. A related
distinct-representative problem occurs in rainbow Hamiltonicity;
Coulson and Perarnau~\cite{14} treat bounded colour classes in Dirac
graphs. Here the local degree and codegree hypotheses guide the
choice of both the cycle and its representative hyperedges.

The first structural ingredient is Theorem~\ref{G:thm:low}.
For any nonnegative integer $T$, the corresponding local Fan and codegree
conditions force Hamiltonicity if a longest Berge path maximizing the
sum of its endpoint degrees has an endpoint of degree at most $T$.
The proof moves the low endpoint through its initial neighbourhood,
then uses $2$-connectivity of the shadow to find two successive
crossing edges. Every exchange explicitly records its released and
inserted representative. No estimate depending on $n$ or on the
formula for $\tau_n$ is used.

For high endpoints, Lemma~\ref{G:lem:endpoint-charging} injectively
assigns incident hyperedges to endpoint-index pairs and a bounded
collection of exceptional slots. It is important that the assignment
is an injection of actual triples, not simply an inequality for their
number. When equality holds, Lemma~\ref{G:lem:one-gap-rigidity}
identifies the exceptional representatives and proves that particular
pairs have codegree one. Condition (b) then makes the opposite
vertices high, and a final rotation contradicts the same endpoint
count. The borderline case $n=10$ requires a separate four-index
replacement to make both endpoint sets prime without losing indices
already secured on the opposite side.

Lemma~\ref{G:lem:uniform-cycle-extension} provides a common
cycle-extension argument for all $n\ge9$. Section~\ref{sec:large-order}
settles $n\ge11$, and Section~\ref{sec:boundary} handles the remaining
critical endpoint configuration at $n=10$. For $n=7,8$, the reduction
in Section~\ref{sec:order7} leads to the exhaustive search specified
and proved complete in Appendix~\ref{app:verification}.

The proof dependencies are deliberately one-way. The longest-path endpoint
package, the representative-exchange ledger, and the charging injection are
established before either high-endpoint theorem. The cycle-extension lemma
is proved independently of those theorems. Whenever a final rotation creates
another longest path, the proof reapplies only the structural lemmas and the
numerical reductions already established for \emph{every} longest path
with high-degree endpoints, not the case analysis being concluded.

\section{Common longest-path tools}\label{sec:common-tools}
We collect the rotation and charging lemmas used throughout the proof.

\begin{remark}[Representative safety]\label{G:rem:rep-safe}
Each switch must preserve distinct key vertices and distinct representative
hyperedges. A hyperedge may be assigned to a new position only after its
previous occurrence has been removed; the nontrivial exchanges below are
checked explicitly.
\end{remark}

\begin{lemma}[Prefix-rotation ledger]\label{G:lem:prefix-rotation-ledger}
Let $Q=u_0g_1u_1\cdots g_mu_m$ be a Berge path of length $m\ge2$,
let $2\le j\le m$, and let $f$ contain $\{u_0,u_j\}$.
Suppose either $f\notin E(Q)$ or $f=g_j$. Then
\[
 u_{j-1}g_{j-1}\cdots g_1u_0\,f\,u_jg_{j+1}\cdots g_mu_m
 \tag{PR}
\]
is a Berge path with the same key-vertex set and the same length as $Q$.
Its representative set is, respectively,
\[
 (E(Q)\setminus\{g_j\})\cup\{f\}
 \quad\text{or}\quad E(Q).
\]
If $j=m$ and $f\notin E(Q)$, the original path $Q$ together
with $f$ is also a Berge cycle of length $m+1$.
\end{lemma}
\begin{proof}
The two portions of \textup{(PR)} use disjoint sets of key vertices:
$u_{j-1},\ldots,u_0$ and $u_j,\ldots,u_m$. They cover $K(Q)$.
All old representatives except $g_j$ occur once, in either their old
or reversed incidence. The inserted incidence $u_0fu_j$ is valid.
If $f$ is outside $E(Q)$ it is new; if $f=g_j$, its old incidence
$u_{j-1}g_ju_j$ has been deleted. In both cases representatives remain
distinct. Finally, for $j=m$ an outside $f$ joins the endpoints of $Q$
without using an old representative.
\end{proof}

\begin{lemma}[Berge distance two and the shadow]\label{G:lem:shadow-distance}
For distinct vertices $u,v$ of a hypergraph $H$,
\[
 d_H(u,v)=2
 \quad\Longleftrightarrow\quad
 \operatorname{dist}_{\partial_2H}(u,v)=2.
\]
\end{lemma}
\begin{proof}
A Berge $2$-path gives a $2$-walk in the shadow. Conversely,
let $u-w-v$ be a shadow $2$-path.  Choose hyperedges
$e\supseteq\{u,w\}$ and $f\supseteq\{w,v\}$.  Since the shadow distance
between $u$ and $v$ is two, $uv$ is not a shadow edge.  Hence $e\ne f$;
otherwise the common hyperedge would contain both $u$ and $v$.  Therefore
$u,e,w,f,v$ is a Berge $2$-path.  The nonadjacency of $u,v$ also excludes a
Berge path of length one.
\end{proof}

\begin{lemma}[Low-degree multiplicity]\label{G:lem:low-multiplicity}
Let $T$ be a nonnegative integer.  Suppose that a simple $3$-uniform
hypergraph $H$ satisfies
\begin{equation}\label{G:cond:C}
 \deg_H(u,v)=1\quad\Longrightarrow\quad
 \min\{\deg_H(u),\deg_H(v)\}\ge T+1.
\end{equation}
If $\deg_H(x)\le T$ and $xy\in E(\partial_2H)$, then
\[
 \deg_H(x,y)\ge2.
\]
\end{lemma}
\begin{proof}
The codegree is positive because $xy$ is a shadow edge.  It cannot equal
one, since then \eqref{G:cond:C} would imply $\deg_H(x)\ge T+1$.
\end{proof}

\begin{lemma}[Cycle-to-path conversion]\label{G:lem:cycle-to-path}
Let $H$ be connected.  If $H$ contains a non-Hamiltonian Berge cycle of
length $\ell$, then $H$ contains a Berge path of length $\ell$.
\end{lemma}
\begin{proof}
Write
\[
 C=v_1e_1v_2e_2\cdots v_\ell e_\ell v_1.
\]
If some representative $e_i$ contains a vertex $x\notin K(C)$, open the
cycle at $v_i$ and use $e_i$ as the final representative.  Otherwise no
representative of $C$ meets $V(H)\setminus K(C)$.  Since $H$ is connected,
some edge $f\notin E(C)$ meets both $K(C)$ and its complement; opening the
cycle at a vertex of $f\cap K(C)$ and using $f$ last again gives a Berge
path of length $\ell$.  Both constructions use distinct key vertices and representatives.
\end{proof}

\begin{corollary}[Longest-path cycle barrier]\label{G:cor:cycle-barrier}
Assume that $H$ is connected and has no Hamiltonian Berge cycle, and let
$P$ be a longest Berge path of length $m$.  Then $H$ has no Berge cycle of
length larger than $m$.
\end{corollary}
\begin{proof}
A longer non-Hamiltonian cycle would, by Lemma~\ref{G:lem:cycle-to-path},
yield a path longer than $P$; Hamiltonian cycles are excluded.
\end{proof}

\begin{lemma}[Representative-safe cycle switches]\label{G:lem:cycle-switches}
Let
\[
 C=v_0e_0v_1e_1\cdots v_{k-1}e_{k-1}v_0
\]
be a Berge cycle, let $x\notin K(C)$, and read indices modulo $k$.
\begin{enumerate}[label=\textup{(\roman*)}]
\item If $f\notin E(C)$ contains $v_i,v_{i+1}$ and both $e_i,e_{i+1}$
contain $x$, then replacing the two consecutive incidences by
\[
 v_i f v_{i+1}e_i x e_{i+1}v_{i+2}
\]
and retaining the complementary arc of $C$ gives a Berge cycle on
$K(C)\cup\{x\}$.
\item If $f\notin E(C)$ contains nonconsecutive rim vertices $v_i,v_j$ and
both $e_i,e_j$ contain $x$, then the two arcs of $C-\{e_i,e_j\}$, together
with $v_{i+1}e_i x e_jv_{j+1}$ and the chord $v_i f v_j$, give a Berge
cycle on $K(C)\cup\{x\}$.
\end{enumerate}
In both switches the representatives are pairwise distinct.
\end{lemma}
\begin{proof}
In (i), the representatives are $E(C)$ together with $f$; the old
incidences of $e_i,e_{i+1}$ are rerouted through $x$.
In (ii) they are $E(C)\setminus\{e_i,e_j\}$ together with $e_i,e_j,f$ in
the new positions.  The representatives remain distinct, and the two retained arcs
partition the old rim vertices.
\end{proof}

\begin{lemma}[Cycle-extension count]\label{G:lem:cycle-count}
Let $H$ be a simple $3$-uniform hypergraph, and let $C=v_0e_0\cdots v_{k-1}e_{k-1}v_0$ be a Berge cycle and
$x\notin K(C)$ with $N_H(x)\subseteq K(C)$.  Assume that no Berge cycle of
length $k+1$ exists, and among the cycles with key-vertex set $K(C)$ choose
$C$ so that the number of its representatives containing $x$ is maximal.
Put
\[
 H_C:=H-E(C),\qquad I:=\{i\in\mathbb Z_k:v_i\in N_{H_C}(x)\},\qquad s:=|I|.
\]
Then $I$ has no two cyclically consecutive indices, $s\ne1$, and
\begin{equation}\label{G:eq:cycle-count}
 \deg_H(x)\le \binom{s}{2}+k-2s.
\end{equation}
Moreover, if $i\in I$, neither $e_{i-1}$ nor $e_i$ contains $x$.
\end{lemma}
\begin{proof}
If $i\in I$, an outside witness for $xv_i$ and an incident cycle
representative containing $x$ would give a $(k+1)$-cycle.  Hence $e_{i-1},e_i$ avoid $x$.  If $i,i+1\in I$, choose
outside witnesses $a,b$ for $xv_i,xv_{i+1}$.  When $a\ne b$, deleting
$e_i$ and inserting $v_i a x b v_{i+1}$ gives a $(k+1)$-cycle; when
$a=b$, the common edge is $\{x,v_i,v_{i+1}\}$ and replacing $e_i$ by it
increases the number of cycle representatives containing $x$.  Hence $I$ is cyclically independent. The case $s=1$ is impossible:
an outside edge through $x$ meets two rim vertices, giving two indices.

Every outside edge through $x$ is $\{x,v_i,v_j\}$ with $i,j\in I$, so
there are at most $\binom{s}{2}$ of them.  For each $i\in I$ the two
representatives $e_{i-1},e_i$ avoid $x$, and these $2s$ representatives are
distinct because $I$ has no consecutive indices.  At most $k-2s$ cycle
representatives can therefore contain $x$, which proves
\eqref{G:eq:cycle-count}.
\end{proof}

\begin{lemma}[An alternating complete-link configuration]\label{G:lem:parity-link}
Let $r\ge3$, put $L=2r+1$, and let $H$ be a simple $3$-graph on
$\{x,v_1,\ldots,v_L\}$. Suppose
\[
 C=v_1e_1v_2e_2\cdots v_Le_Lv_1
\]
is a Berge cycle, $e_1=\{x,v_1,v_2\}$, and, with
$O=\{v_3,v_5,\ldots,v_L\}$,
\[
 E_H(x)=\{e_1\}\cup\{\{x,u,v\}:u,v\in O,\ u\ne v\}.
\]
If $H$ has no Hamiltonian Berge cycle, then
$N_H(v_q)\subseteq O$ for every even $q\in\{4,6,\ldots,2r\}$.
\end{lemma}
\begin{proof}
Every prescribed edge through $x$ other than $e_1$ lies outside $E(C)$,
since its two rim vertices are nonconsecutive. Each vertex of $O$ lies in
$r-1\ge2$ such edges. Thus outside witnesses at two distinct vertices
of $O$ can be chosen distinct; their candidate families have just one
common edge. No such witness contains an interior even vertex.

We use two switches. If $f$ contains nonconsecutive
$v_i,v_j$, deleting $e_i,e_j$ and choosing outside witnesses $a,b$ at
$v_{i+1},v_{j+1}$ gives
\[
 v_{i+1}a x b v_{j+1}\,C\,v_i f v_j\,C^{-1}\,v_{i+1}.
 \tag{$+$}
\]
Deleting $e_{i-1},e_{j-1}$ and using witnesses at the predecessors gives
\[
 v_{i-1}a x b v_{j-1}\,C^{-1}\,v_i f v_j\,C\,v_{i-1}.
 \tag{$-$}
\]
Cycle indices in these two formulas are read modulo $L$. In each formula
the two retained arcs partition the rim vertices. The exact representative ledgers are
\[
\begin{array}{ll}
(+):&(E(C)\setminus\{e_i,e_j\})\cup\{a,b,f\},\\[1mm]
(-):&(E(C)\setminus\{e_{i-1},e_{j-1}\})\cup\{a,b,f\}.
\end{array}
\]
The constructions are safe when $f$ is outside $C$ or is one of the
two representatives deleted in the chosen switch, and $a,b,f$ are
distinct. In particular, using a deleted representative as the chord
does not use a second copy of that hyperedge.

Suppose first that $f$ contains two distinct interior even vertices
$v_q,v_s$. If $f\notin E(C)$ or $f\in\{e_q,e_s\}$, use $(+)$.
Otherwise $3$-uniformity forces
$f\in\{e_{q-1},e_{s-1}\}$, and use $(-)$. All attachment vertices
belong to $O$. The complete link supplies distinct $a,b$, both different from $f$
because $f$ meets an interior even vertex. Thus either construction is Hamiltonian, a contradiction.

Next let $f$ contain $v_1,v_q$, with $q$ interior and even. A cycle
representative containing this pair must belong to
$\{e_L,e_{q-1},e_q\}$, since $e_1$ is prescribed. If $f$ is outside
$C$ or belongs to $\{e_L,e_{q-1}\}$, use $(-)$ with $i=1,j=q$;
the attachment vertices $v_L,v_{q-1}$ lie in $O$. For the remaining
case $f=e_q=\{v_q,v_{q+1},v_1\}$, choose an outside witness $a$ at
$v_{q+1}$ and use
\[
 v_1e_Lv_Le_{L-1}\cdots e_{q+1}v_{q+1}a x e_1v_2e_2\cdots
 e_{q-1}v_q f v_1.
\]
The descending segment is empty when $q=L-1$. Here $f$ is removed from
its old position, $a$ lies outside $E(C)$, and all other representatives
are distinct cycle edges. The representatives of this exceptional cycle are exactly
$E(C)\cup\{a\}$; notably $f=e_q$ appears only in its final chord
position. Hence $v_1$ is not adjacent to $v_q$.

Finally, a cycle representative containing $v_2,v_q$ belongs to
$\{e_2,e_{q-1},e_q\}$. An outside $f$ or
$f\in\{e_2,e_q\}$ permits $(+)$ with $i=2,j=q$, using the distinct
attachment vertices $v_3,v_{q+1}\in O$. If instead
$f=e_{q-1}=\{v_{q-1},v_q,v_2\}$, choose an outside witness $a$ at
$v_{q-1}$ and use
\[
 v_1e_1x a v_{q-1}e_{q-2}\cdots e_2v_2 f v_qe_qv_{q+1}\cdots e_Lv_1.
\]
Again $f$ replaces its old incidence, $a$ is new, and all representatives
are distinct. This second exceptional cycle also has representative set
$E(C)\cup\{a\}$; here $f=e_{q-1}$ occurs only on $v_2v_q$.
Thus $v_q$ has no neighbour among $v_1,v_2$, the other
interior even vertices, or $x$ (the last exclusion follows from the
prescribed link). All remaining vertices belong to $O$.
\end{proof}

\begin{lemma}[Few exceptional cycle representatives]\label{G:lem:few-exceptions}
Let $T$ be a nonnegative integer and let $H$ be a simple $3$-graph satisfying
\[
 d_H(u,v)=2\quad\Longrightarrow\quad
 \max\{\deg_H(u),\deg_H(v)\}\ge T+1.
\]
Let $C=v_1e_1\cdots v_ke_kv_1$ be a Berge cycle with $k\ge7$, and let
$x\notin K(C)$ satisfy $E_H(x)\subseteq E(C)$. Put
\[
 b=|\{i:x\notin e_i\}|,
 \qquad h=|V(H)\setminus(K(C)\cup\{x\})|.
\]
If
\[
 k>3b\qquad\text{and}\qquad
 2+b+\binom{b+h}{2}\le T,
\]
then $H$ contains a Berge cycle of length $k+1$.
\end{lemma}
\begin{proof}
Assume that no such cycle exists. Call $i$ good if $x\in e_i$.
An edge outside $E(C)$ cannot contain two good-index rim vertices:
Lemma~\ref{G:lem:cycle-switches}(i) or (ii) would insert $x$.
It avoids $x$ by hypothesis. A good-index vertex therefore belongs
to at most two incident cycle representatives, $b$ other cycle
representatives, and $\binom{b+h}{2}$ outside edges. Indeed, the other vertices of an outside edge through it must
be chosen among the $b$ bad-index rim vertices and the $h$ remaining
vertices. Its degree is therefore at most $T$.

Read rim indices modulo $k$. Among the $k$ pairs $\{v_i,v_{i+2}\}$,
at least $k-2b$ have both
indices good. A good representative contains no such pair, and neither
does an outside edge. Each bad representative contains its two consecutive
rim vertices and only one further vertex; since $k\ge7$, it contains
at most one pair of cyclic distance two. Hence at least $k-3b>0$ of
the good-good pairs are nonadjacent. For any such pair,
$v_ie_iv_{i+1}e_{i+1}v_{i+2}$ is a Berge path with two distinct
representatives. Its endpoints have degree at most $T$, contradicting
the distance-two hypothesis.
\end{proof}

\subsection{Order-independent endpoint package}
\begin{lemma}\label{G:lem:endpoint-package}
Let $H$ be a connected simple $3$-uniform hypergraph with no Hamiltonian
Berge cycle, and let
\[
 P=v_0e_1v_1\cdots e_mv_m
\]
be a longest Berge path.  Set $H'=(V(H),E(H)\setminus E(P))$.  All indices
below lie in $\{0,1,\ldots,m-1\}$.  Define
\[
\begin{aligned}
I_1 &:= \{\, i : v_{i+1} \in N_{H'}(v_0) \text{ or } e_{i+1} = \{v_0, v_i, v_{i+1}\} \},\\
I_1' &:= \{\, i : v_{i+1} \in N_{H'}(v_0) \},\\
I_2 &:= \{\, i : v_i \in N_{H'}(v_m) \text{ or } e_{i+1} = \{v_m, v_i, v_{i+1}\} \},\\
I_2' &:= \{\, i : v_i \in N_{H'}(v_m) \},\\
I_3 &:= \{\, i : i \in I_1' \text{ and } i-1 \in I_2 \},\\
I_4 &:= \{\, i : i+1 \in I_1 \text{ and } i \in I_2' \}.
\end{aligned}
\]
Then:
\begin{enumerate}
\item \(I_1 \cap I_2 = \emptyset\). Consequently, \(|I_1| + |I_2| \le m\).
\item
\[
\deg(v_0) \le \binom{|I_1'|}{2} + |I_1| + 1.
\]
\item
\[
\deg(v_m) \le \binom{|I_2'|}{2} + |I_2| + 1.
\]
\item For every \(i \in I_3\) (i.e., \(i \in I_1'\) and \(i-1 \in I_2\)), if \(e_{i+1} = \{v_0, v_i, v_{i+1}\}\), then \(\deg(v_i, v_{i+1}) = 1\) and \(N(v_i) \subseteq K(P)\).
\item For every \(i \in I_4\) (i.e., \(i \in I_2'\) and \(i+1 \in I_1\)), if \(e_{i+1} = \{v_i, v_{i+1}, v_m\}\), then
\[
\deg(v_i,v_{i+1})=1
\qquad\text{and}\qquad
N_H(v_{i+1})\subseteq K(P).
\]
\item If \(0 \in I_1'\), then \(N_H(v_0) \subseteq K(P)\); if \(m-1 \in I_2'\), then \(N_H(v_m) \subseteq K(P)\).
\end{enumerate}

\end{lemma}
\begin{proof}
Throughout we use Corollary~\ref{G:cor:cycle-barrier}; a path segment
of length zero consists of its endpoint.

\emph{(1).}
If $i\in I_1\cap I_2$, choose witnesses $a,b$ for $v_0v_{i+1}$
and $v_mv_i$, respectively, using outside representatives for prime
indices and $e_{i+1}$ otherwise.  For $1\le i\le m-2$, equality
$a=b$ would put four distinct vertices in one triple.  At $i=0$
or $m-1$, equality would give an outside representative containing
both endpoints and close $P$.  Thus $a\ne b$, and
\[
 v_0e_1\cdots e_iv_i\,b\,v_me_m\cdots
 e_{i+2}v_{i+1}\,a\,v_0
\]
is a Berge cycle of length $m+1$.  Its new representatives are outside
$E(P)$ or equal to the omitted edge $e_{i+1}$, so none is repeated.
Hence $I_1\cap I_2=\varnothing$.

\emph{(2), (3), and (6).}
An edge through $v_0$ meeting $V(H)\setminus K(P)$ must be $e_1$:
every other such edge lies outside $E(P)$ and extends $P$.
Every outside edge through $v_0$ is therefore determined by a pair
in $I_1'$, whereas each path representative other than $e_1$
containing $v_0$ has its index in $I_1$.  This gives (2), and path
reversal gives (3).  If $0\in I_1'$, choose an outside witness $a$
for $v_0v_1$.  Were $e_1=\{v_0,v_1,u\}$ with $u\notin K(P)$, then
$u e_1v_0 a v_1e_2\cdots e_mv_m$ would extend $P$.
Thus $N_H(v_0)\subseteq K(P)$; reversal proves the other half of (6).

\emph{(4).}
Fix $i\in I_3$ and suppose $e_{i+1}=\{v_0,v_i,v_{i+1}\}$.
Choose $a\in E(H')$ containing $v_0,v_{i+1}$ and a witness $b$
for $i-1\in I_2$.  Thus $b$ is outside $E(P)$ or
$b=e_i=\{v_m,v_{i-1},v_i\}$.  If $i=m-1$, then $a$ joins the endpoints $v_0,v_m$ outside
$E(P)$ and closes $P$, so $1\le i\le m-2$.  Also $a\ne b$:
for $i\ge2$ equality would force the four distinct vertices
$v_0,v_{i+1},v_m,v_{i-1}$ into one triple; for $i=1$ it would
produce an outside edge containing $v_0,v_m$ and close $P$.

Suppose $f\ne e_{i+1}$ contains $v_i,v_{i+1}$.  By uniformity,
$f\notin E(P)$ or $f\in\{e_i,e_{i+2}\}$.  If $f\ne e_{i+2}$, use
\[
 v_0e_1\cdots e_{i-1}v_{i-1}\,b\,v_me_m\cdots
 e_{i+2}v_{i+1}\,f\,v_i e_{i+1}v_0.
\]
Here $f\ne b$: if both are outside, equality would require the
four distinct vertices $v_i,v_{i+1},v_m,v_{i-1}$;
if $f=b=e_i$, it would force $v_{i+1}=v_m$, or $i=m-1$.
Moreover, if $f=a$, then $f=\{v_0,v_i,v_{i+1}\}=e_{i+1}$,
contrary to the choice of $f$.  The displayed cycle does not use
$a$; the last observation also makes the witness separation explicit.
All representatives used in this cycle are distinct.  If $f=e_{i+2}$, use
\[
 v_0e_1\cdots e_{i-1}v_{i-1}\,b\,v_me_m\cdots
 e_{i+3}v_{i+2}\,e_{i+2}v_i e_{i+1}v_{i+1}\,a\,v_0.
\]
In this construction the old incidence of $e_i$ is deleted;
if $b=e_i$, that representative is used only in the new position
between $v_{i-1}$ and $v_m$.  The outside witness $a$ is not a path
representative, and $b$ is either outside $E(P)$ or the released $e_i$.
The two witnesses are distinct, so all cycle representatives are
pairwise distinct, including the boundary case $i=m-2$ where the
descending segment is empty. Both constructions give a forbidden
$(m+1)$-cycle.  Hence $\deg_H(v_i,v_{i+1})=1$.

Finally, suppose $u\notin K(P)$ and $h\supseteq\{u,v_i\}$.
Uniformity and the displayed form of $e_{i+1}$ imply that
$h\notin E(P)$ or $h=e_i$.  In the first case,
\[
 u h v_i e_i v_{i-1}\cdots e_1v_0
 e_{i+1}v_{i+1}e_{i+2}\cdots e_mv_m
\]
extends $P$.  In the second case $b$ must be outside $E(P)$: the edge
$e_i$ contains the non-key vertex $u$ and hence cannot equal the
special edge $\{v_m,v_{i-1},v_i\}$ used as a path-only witness
for $i-1\in I_2$. Therefore
\[
 u e_i v_i e_{i+1}v_{i+1}\cdots e_mv_m
 b v_{i-1}e_{i-1}\cdots e_1v_0
\]
does so.  Thus $N_H(v_i)\subseteq K(P)$.

\emph{(5).}
Reverse $P$ and apply (4) with $j=m-i-1$.
The conditions $i\in I_2'$ and $i+1\in I_1$ become
$j\in I_1'$ and $j-1\in I_2$, respectively, and the central-edge
hypothesis becomes that of (4).
\end{proof}

\begin{lemma}[First-edge crossing principle]\label{G:lem:first-edge-crossing}
Let $P=v_0e_1v_1\cdots e_mv_m$ be a Berge path of length
$m\ge2$, and let $I_2$ and $I_1'$ be its endpoint index sets.
Suppose $0\in I_1'$ and
\[
 e_1=\{v_0,v_1,v_{j+1}\},\qquad 1\le j\le m-1.
\]
If $j\in I_2$, then $H$ contains a Berge cycle of length $m+1$.
\end{lemma}
\begin{proof}
Choose an edge $a\notin E(P)$ containing $v_0,v_1$, and a witness
$b$ for $j\in I_2$.  Thus either $b\notin E(P)$ contains $v_m,v_j$,
or $b=e_{j+1}=\{v_m,v_j,v_{j+1}\}$.
If $a=b$, then $b$ must be outside $E(P)$; since a triple contains
$v_0,v_1,v_m,v_j$, with $1\le j<m$, necessarily $j=1$.
Consequently $a$ contains both $v_0,v_m$ and closes $P$.
If $a\ne b$, use
\[
 v_0\,a\,v_1e_2\cdots e_jv_j\,b\,v_m
 e_m\cdots e_{j+2}v_{j+1}\,e_1\,v_0.
\]
The retained path representatives are $E(P)\setminus\{e_{j+1}\}$,
while $a$ lies outside $E(P)$ and $b$ is outside $E(P)$ or equals
$e_{j+1}$.  The three newly used incidences have pairwise distinct
representatives, and the two arcs visit all $m+1$ key vertices once.
\end{proof}

\begin{lemma}[First-representative exchange]
\label{G:lem:first-edge-ledger}
Let $P=v_0e_1v_1\cdots e_mv_m$ be a Berge path and suppose
$a=\{v_0,v_1,v_{j+1}\}\in E(H)\setminus E(P)$ for some
$1\le j\le m-1$.  Replace $e_1$ by $a$ to obtain $P^*$.
Then $P^*$ is a Berge path with the same key order and
\[
 E(P^*)=(E(P)\setminus\{e_1\})\cup\{a\}.
\]
Every old left outside-witness index except possibly $0,j$ survives.
If also $e_1=\{v_0,v_1,v_{\ell+1}\}$ for an index
$1\le\ell\le m-1$, then $0,\ell\in I_1'(P^*)$.
If $v_m\notin a$, then
\[
 I_2'(P)\subseteq I_2'(P^*),\qquad
 I_2(P)\subseteq I_2(P^*).
\]
If instead $v_m\in a$, then $P$ immediately closes to a Berge
cycle of length $m+1$.
\end{lemma}
\begin{proof}
The displayed representative-set identity follows because $a$ is
outside $E(P)$, so all representatives in $P^*$ are distinct.
The outside-edge sets satisfy the exact identity
\[
 E(H)\setminus E(P^*)=
 (E(H)\setminus E(P))\setminus\{a\}\,\cup\,\{e_1\}.
\]
Since $a$ has non-$v_0$ vertices $v_1,v_{j+1}$, only left
outside witnesses at indices $0,j$ can be lost. If $e_1$ has
the stated form, the released representative witnesses $0,\ell$.
In particular, preservation of a prime index $i\notin\{0,j\}$
is literal preservation of an outside triple; no new representative
is being assumed.
If $v_m\notin a$, no old outside edge incident with $v_m$ is
removed, yielding the first right-side inclusion.  Every old
path-only right witness at index $i\ge1$ is unchanged.
The sole remaining case is a path-only right witness at index
$i=0$, which requires $e_1=\{v_m,v_0,v_1\}$; after the exchange
this old $e_1$ lies outside $P^*$ and witnesses the same index
through the pair $v_m,v_0$.  Thus all old right indices survive.
When $v_m\in a$, the representative $a$ directly closes $P$
without repeating a path representative.
\end{proof}

\begin{lemma}[High-endpoint structural transfer]\label{G:lem:high-structure}
Let $T\ge0$, let $H$ be a connected simple $3$-uniform hypergraph with no
Hamiltonian Berge cycle, and let
$P=v_0e_1v_1\cdots e_mv_m$ be a longest Berge path with
$\deg_H(v_0),\deg_H(v_m)\ge T+1$.  Assume
\[
 \deg_H(u,v)=1\Longrightarrow
 \min\{\deg_H(u),\deg_H(v)\}\ge T+1,
 \tag{$C_T$}
\]
and assume the following cycle-extension property at length $m$: whenever
$C$ is a Berge cycle of length $m$ and $x\notin K(C)$ satisfies
$N_H(x)\subseteq K(C)$ and $\deg_H(x)\ge T+1$, the hypergraph contains a
Berge cycle longer than $m$.

Then, for the endpoint sets of Lemma~\ref{G:lem:endpoint-package},
\begin{enumerate}[label=\textup{(\roman*)}]
\item $I_1\cap I_2=\varnothing$;
\item $0\notin I_2$ and $m-1\notin I_1$;
\item if $i\in I_1'$ and $i-1\in I_2$, then
$e_{i+1}\ne\{v_0,v_i,v_{i+1}\}$;
\item if $i\in I_2'$ and $i+1\in I_1$, then
$e_{i+1}\ne\{v_i,v_{i+1},v_m\}$;
\item
\[
 \deg_H(v_0)\le\binom{|I_1'|}{2}+|I_1|+1,
 \qquad
 \deg_H(v_m)\le\binom{|I_2'|}{2}+|I_2|+1.
\]
\end{enumerate}
\end{lemma}
\begin{proof}
Parts (i) and (v) are Lemma~\ref{G:lem:endpoint-package}(1)--(3).
For (ii), if $0\in I_2'$, an outside edge closes $P$.  Otherwise
$e_1=\{v_m,v_0,v_1\}$, and
\[
 v_1e_2v_2\cdots e_mv_me_1v_1
\]
is a Berge cycle of length $m$ on $K(P)\setminus\{v_0\}$.  Maximality
of \(P\) forbids an edge outside \(E(P)\) through \(v_0\) from meeting a
vertex outside \(K(P)\).  Among the path representatives through \(v_0\),
only \(e_1\) could have such a third vertex; here
\(e_1=\{v_m,v_0,v_1\}\).  Hence \(N_H(v_0)\subseteq K(P)\), and the
extension property yields a cycle longer than \(m\), contrary to
Corollary~\ref{G:cor:cycle-barrier}.  Reversal gives \(m-1\notin I_1\).

For (iii), suppose the forbidden equality holds.  Lemma~\ref{G:lem:endpoint-package}(4)
gives
$\deg_H(v_i,v_{i+1})=1$ and $N_H(v_i)\subseteq K(P)$; hence $(C_T)$ gives
$\deg_H(v_i)\ge T+1$.  With the witnesses $a,b$ from that lemma,
\[
 v_0e_1\cdots e_{i-1}v_{i-1}\,b\,v_m e_m\cdots
 e_{i+2}v_{i+1}\,a\,v_0
\]
is a representative-safe Berge cycle of length $m$ on
$K(P)\setminus\{v_i\}$, which contains every neighbor of $v_i$.  The
extension property again contradicts Corollary~\ref{G:cor:cycle-barrier}.
Part (iv) is the reversed
form of (iii).
\end{proof}

\begin{lemma}[Full endpoint partition]\label{G:lem:full-partition}
Assume the hypotheses of Lemma~\ref{G:lem:high-structure}. If
\[
 I_1\cup I_2=\{0,\ldots,m-1\},\qquad I_1=I_1',
\]
then
\[
 \deg_H(v_0)\le\binom{|I_1|}{2}.
\]
The symmetric statement holds at $v_m$ when $I_2=I_2'$.
\end{lemma}
\begin{proof}
The boundary exclusion $0\notin I_2$ gives $0\in I_1'$.
Lemma~\ref{G:lem:endpoint-package}(6) yields
$N_H(v_0)\subseteq K(P)$, so
$e_1=\{v_0,v_1,v_{j+1}\}$ for some $1\le j\le m-1$.
If $j\in I_2$, Lemma~\ref{G:lem:first-edge-crossing} produces an
$(m+1)$-cycle, contradicting Corollary~\ref{G:cor:cycle-barrier}.
The partition therefore gives $j\in I_1$.

Assign each outside edge through $v_0$ to the pair of predecessor
indices of its other two vertices in $I_1$. For a path representative
$e_{i+1}$ through $v_0$, with $i\ge1$, we have $i\in I_1=I_1'$.
Lemma~\ref{G:lem:high-structure}(iii) excludes $i-1\in I_2$;
the partition forces $i-1\in I_1$, so assign $e_{i+1}$ to
$\{i-1,i\}$. Finally assign $e_1$ to $\{0,j\}$.
Each pair determines its full triple through $v_0$, and simplicity
therefore makes this assignment injective. Reversing the path gives
the right-endpoint statement.
\end{proof}

\begin{lemma}[Endpoint charging]\label{G:lem:endpoint-charging}
Let $H$ be a simple $3$-uniform hypergraph, and let $P=v_0e_1v_1\cdots e_mv_m$ be a longest Berge path, with endpoint sets
$I_1,I_1',I_2,I_2'$ as in Lemma~\ref{G:lem:endpoint-package}.  Assume
$I_1\cap I_2=\varnothing$ and the following crossing restrictions:
\begin{align*}
 i\in I_1',\ i\ge1,\ v_0\in e_{i+1}&\Longrightarrow i-1\notin I_2,\\
 i\in I_2',\ i\le m-2,\ v_m\in e_{i+1}&\Longrightarrow i+1\notin I_1.
\end{align*}
Put
\[
 a=|I_1|,\quad a'=|I_1'|,\qquad b=|I_2|,\quad b'=|I_2'|,
 \qquad M=\{0,\ldots,m-1\}\setminus(I_1\cup I_2),\quad q=|M|.
\]
Then
\begin{equation}\label{G:eq:endpoint-charging}
 \deg_H(v_0)\le\binom{a'}2+2(a-a')+q+1,
 \qquad
 \deg_H(v_m)\le\binom{b'}2+2(b-b')+q+1.
\end{equation}
For the left endpoint put
\[
\begin{aligned}
 \mathcal P_0&=\binom{I_1'}2, &
 \mathcal D_0&=\{D_i:i\in I_1\setminus I_1'\}, &
 \mathcal R_0&=\{R_i:i\in I_1\setminus I_1'\},\\
 \mathcal M_0&=\{M_i:i\in M\}, &
 \mathcal S_0&=\{\star_0\}.&&
\end{aligned}
\]
There is an explicit injection
\[
 \Phi_0:E_H(v_0)\hookrightarrow
 \mathcal P_0\dot\cup\mathcal D_0\dot\cup\mathcal R_0
 \dot\cup\mathcal M_0\dot\cup\mathcal S_0.
\]
Here a pair-slot $\{i,j\}\in\mathcal P_0$ refers to the
\emph{actual triple} $\{v_0,v_{i+1},v_{j+1}\}$;
the other slots record individual path representatives.
The right endpoint has the reverse-path analogue.  Hence equality in
the corresponding bound makes the injection bijective and every target slot
is occupied.
\end{lemma}
\begin{proof}
We prove the left inequality.  Every edge outside $E(P)$ through $v_0$ has
its other two vertices in $K(P)$ and is therefore uniquely charged, by
simplicity, to the corresponding pair in $\binom{I_1'}2$.

For a path representative $e_{i+1}$ containing $v_0$, charge $e_1$ to
$\star_0$.  If $i\in I_1\setminus I_1'$, charge it to $D_i$.  If
$i\in I_1'$ and $i\ge1$, then the crossing restriction gives
$i-1\notin I_2$.  Charge $e_{i+1}$ to $\{i-1,i\}$ when
$i-1\in I_1'$, to $R_{i-1}$ when $i-1\in I_1\setminus I_1'$, and to $M_{i-1}$ otherwise.  In the pair-slot case,
$e_{i+1}=\{v_0,v_i,v_{i+1}\}$, precisely the triple
represented by $\{i-1,i\}$. Thus an outside edge cannot
occupy the same pair-slot without coinciding with a path edge,
contrary to its being outside $E(P)$. Different pair-slots cannot
represent the same triple, since their two non-$v_0$ vertices
uniquely determine the index pair. The labels of $D_i,R_i,M_i$
uniquely determine the relevant path edge $e_{i+1}$, while
$\star_0$ contains only $e_1$; the slot classes are disjoint.
This proves the injection and the left bound. Path reversal gives the
right bound. Equality makes the finite injection bijective.
\end{proof}

\begin{lemma}[One-gap charging rigidity]
\label{G:lem:one-gap-rigidity}
Assume the hypotheses and crossing restrictions of
Lemma~\ref{G:lem:endpoint-charging}.  Suppose that
\[
 I_1=I_1'=S,\quad 0\in S,\qquad
 \{0,\ldots,m-1\}\setminus(I_1\cup I_2)=\{k\},
 \quad 0<k<m-1,
\]
that $k+1\in S$, and that the first path representative is
\emph{pair-chargeable}, in the explicit sense
\[
 e_1=\{v_0,v_1,v_{j+1}\}\quad\text{for some }j\in S\setminus\{0\}.
\]
Then there is an injection
\[
 E_H(v_0)\hookrightarrow \binom S2\,\dot\cup\,\{\mu_k\},
\]
where $\mu_k$ can be occupied only by
$e_{k+2}=\{v_0,v_{k+1},v_{k+2}\}$.  Consequently
\[
 \deg_H(v_0)\le\binom{|S|}{2}+1.
\]
If $\deg_H(v_0)\ge\binom{|S|}{2}+1$, equality holds, every
$S$-pair slot is occupied, and $e_{k+2}$ has the displayed form.
In this equality case,
\[
 \deg_H(v_0,v_{k+1})=1.
\]
The reversed-path statement is valid with $I_2=I_2'$ and
$k-1\in I_2$; its exceptional representative is
$e_k=\{v_m,v_{k-1},v_k\}$ and the unique codegree-one pair is
$\{v_m,v_k\}$.
\end{lemma}
\begin{proof}
An edge outside $E(P)$ containing $v_0$ has two other vertices
$v_{i+1},v_{j+1}$ with $i,j\in S$, and occupies the pair-slot
$\{i,j\}$.  The first path representative occupies $\{0,j\}$
by the pair-chargeability hypothesis.  Any other path representative
containing $v_0$ is $e_{i+1}=\{v_0,v_i,v_{i+1}\}$ with
$i\in S$ and $i\ge1$.  By the crossing restriction,
$i-1\notin I_2$.  As there is exactly one missing index,
$i-1\in S\cup\{k\}$.  In the former case charge $e_{i+1}$
to $\{i-1,i\}\in\binom S2$.  In the latter case necessarily
$i=k+1$, and the only such path representative is $e_{k+2}$;
charge it to $\mu_k$.  An $S$-pair determines the full triple
through $v_0$, so distinct edges cannot collide in a pair-slot;
$\mu_k$ is formally separate.  This proves injectivity.

If the degree reaches $\binom{|S|}2+1$, each target slot must be
occupied. In particular $e_{k+2}$ has the specified form.
Since $k\notin I_1$, no outside hyperedge contains
$\{v_0,v_{k+1}\}$.  Of the path edges only
$e_1,e_{k+1},e_{k+2}$ could contain that pair: every other
path edge already contains two key vertices other than $v_0$
and $v_{k+1}$.  The first is excluded by pair-chargeability
($j\in S$, $k\notin S$); the second is excluded because
$v_0\in e_{k+1}$ would put $k\in I_1$ via its special-edge
alternative.  Thus $e_{k+2}$ is its unique representative.
Reverse the vertex and edge orders to get the final assertion.
\end{proof}

\begin{theorem}[Order-independent low-degree endpoint theorem]\label{G:thm:low}
Let $T$ be a nonnegative integer and let $H$ be a simple
$3$-uniform hypergraph whose $2$-shadow is $2$-connected.  Assume that for
all distinct vertices $u,v$,
\begin{align}
 d_H(u,v)=2
 &\quad\Longrightarrow\quad
 \max\{\deg_H(u),\deg_H(v)\}\ge T+1,
 \tag{F$_T$}\label{G:cond:F}\\
 \deg_H(u,v)=1
 &\quad\Longrightarrow\quad
 \min\{\deg_H(u),\deg_H(v)\}\ge T+1.
 \tag{C$_T$}
\end{align}
Among all longest Berge paths choose
\[
 P=v_0e_1v_1\cdots e_mv_m
\]
maximizing $\deg_H(v_0)+\deg_H(v_m)$.  If
\[
 \min\{\deg_H(v_0),\deg_H(v_m)\}\le T,
\]
then $H$ contains a Hamiltonian Berge cycle.
\end{theorem}
\begin{proof}
Suppose that $H$ has no Hamiltonian Berge cycle, and assume
$\deg_H(v_0)\le T$. The $2$-connected shadow makes $H$ connected,
so Corollary~\ref{G:cor:cycle-barrier} excludes every Berge cycle of
length $m+1$.  Every vertex has a shadow neighbour.  In particular,
if $x$ is low and $xy$ is a shadow edge, condition $(C_T)$ forces
$\deg_H(x,y)\ge2$, and hence $\deg_H(x)\ge2$.
Thus the $2$-edge-connectivity hypothesis sometimes included in the
low-endpoint reduction is unnecessary.
Call a vertex \emph{low}
if its degree is at most $T$.  By $(F_T)$ and
Lemma~\ref{G:lem:shadow-distance},
\begin{equation}\label{G:eq:low-common-neighbour}
 \begin{gathered}
 \text{Two distinct low vertices with a common shadow neighbour}\\
 \text{are adjacent in the shadow.}
 \end{gathered}
\end{equation}

If a longest path $Q$ has key set $K(P)$ and a low endpoint
$x$, then $N_H(x)\subseteq K(P)$. Indeed, for any non-key neighbour
$u$, Lemma~\ref{G:lem:low-multiplicity} supplies at least two
representatives of $xu$.  At most the first edge of $Q$, oriented
from $x$, can contain both vertices; an outside representative
therefore extends $Q$.

Let $k=\max\{q:v_0v_q\in E(\partial_2H)\}$.  Since
$\deg_H(v_0)\ge2$, simplicity gives at least three shadow neighbours,
and the preceding observation gives $k\ge3$.  We also have $k<m$.
Otherwise no outside edge contains $v_0,v_m$, and multiplicity
forces
\[
 e_1=\{v_0,v_1,v_m\},\qquad
 e_m=\{v_0,v_{m-1},v_m\}.
\]
An outside edge containing $v_0,v_1$, with $e_1$ as a closing
representative, would yield an $(m+1)$-cycle.  Thus multiplicity
forces $e_2=\{v_0,v_1,v_2\}$.  The longest path
\[
 v_1e_1v_0e_2v_2e_3\cdots e_mv_m
\]
makes $v_1$ low by endpoint maximality.  But $v_1v_m$ has no outside
representative, since one would close this path, and among its
path representatives only $e_1$ contains that pair.  This
contradicts low-degree multiplicity.

\paragraph{Rotations within the initial block.}
We prove backwards that, for $1\le r<k$, $v_0v_r$ is a shadow edge
and a longest path $Q_r$ has endpoints $v_r,v_m$, key set $K(P)$,
and the original terminal segment
$v_ke_{k+1}\cdots e_mv_m$.  Since $P$ was chosen with maximum
endpoint-degree sum among \emph{all} longest paths, and $Q_r$ has
the same terminal endpoint $v_m$, we have
\[
 \deg_H(v_r)+\deg_H(v_m)
 \le \deg_H(v_0)+\deg_H(v_m),\qquad
 \deg_H(v_r)\le\deg_H(v_0)\le T.
\]
Suppose $2\le i\le k$ and $v_0v_i$ is a shadow edge.  An outside
representative $f$ for this pair gives the rotation
\[
 v_{i-1}e_{i-1}\cdots e_1v_0 f v_i e_{i+1}\cdots e_mv_m.
 \tag{R1}
\]
The same construction works with $f=e_i$ when $v_0\in e_i$,
since its old occurrence is omitted.  If neither case applies,
multiplicity and uniformity force $i<k$ and
\begin{equation}\label{G:eq:block-exception}
 E(v_0,v_i)=\{e_1,e_{i+1}\},\qquad
 e_1=\{v_0,v_1,v_i\},\qquad
 e_{i+1}=\{v_0,v_i,v_{i+1}\}.
\end{equation}
Successive applications of multiplicity produce a block
$e_{q+1}=\{v_0,v_q,v_{q+1}\}$ for $i\le q\le s$, followed
by an outside representative $f$ of $v_0v_{s+1}$, with $s<k$.
Such an $f$ must occur: otherwise the block would reach
$e_{k+1}=\{v_0,v_k,v_{k+1}\}$, contrary to the choice of $k$.
We then use
\[
 v_{i-1}e_{i-1}\cdots e_2v_1 e_1v_iP v_s
 e_{s+1}v_0 f v_{s+1}e_{s+2}\cdots e_mv_m.
 \tag{R2}
\]
A segment of length zero consists of its single endpoint;
in particular, the initial segment is just $v_1$ when $i=2$,
and $v_iP v_s$ is just $v_i$ when $s=i$.
The representative set in (R2) is exactly
\[
 E(Q_{i-1})=(E(P)\setminus\{e_i\})\cup\{f\}.
\]
The representatives remain distinct: $f\notin E(P)$,
$e_i$ is omitted, and every remaining old edge appears once.
In particular, the vertices in the initial reversed block and the
unchanged tail are disjoint. The incidence changes in (R2) are
particularly important: $e_1$ is rerouted from $v_0v_1$ to $v_1v_i$,
$e_{s+1}$ is rerouted from $v_sv_{s+1}$ to $v_sv_0$, and $f$ is used on
$v_0v_{s+1}$; the old incidence of $e_i$ is deleted. These are changes of
position, not additional copies of any representative. The full record is
\[
\begin{array}{c|c|c|c}
 &\text{deleted old incidence}&\text{inserted incidence}&\text{representative set}\\ \hline
 (\mathrm{R1}),\ f\notin E(P)&e_i:v_{i-1}v_i&f:v_0v_i&(E(P)\setminus\{e_i\})\cup\{f\}\\
 (\mathrm{R1}),\ f=e_i&e_i:v_{i-1}v_i&e_i:v_0v_i&E(P)\\
 (\mathrm{R2})&e_i:v_{i-1}v_i&f:v_0v_{s+1}&(E(P)\setminus\{e_i\})\cup\{f\}
\end{array}
\]
In (R2), the reroutings of $e_1$ and $e_{s+1}$ are already accounted for
by the displayed vertex sequence; neither edge is duplicated. Both
rotations retain the required terminal segment and make $v_{i-1}$ low.
The vertices $v_0,v_{i-1}$ are adjacent when $f=e_i$;
otherwise they have a common shadow neighbour ($v_i$ in
(R1), or $v_i$ via $e_1,e_i$ in (R2)).
Equation~\eqref{G:eq:low-common-neighbour} gives their adjacency.
This completes the induction. The construction also shows that
the first edge of $Q_r$ joins $v_r$ to $v_{r-1}$ for $r\ge2$.
For $r=1$, it either joins $v_1$ to $v_0$ or is the exceptional
edge $e_1=\{v_0,v_1,v_2\}$ in (R2).

\paragraph{Crossing the block.}
Every vertex of $A=\{v_0,\ldots,v_{k-1}\}$ is low and is an endpoint
of a longest path with key set $K(P)$; hence none has a neighbour
outside $K(P)$.  The $2$-connectivity of the shadow supplies a crossing edge
$v_iv_j$ with $1\le i<k<j\le m$. Indeed, if no such edge
existed, deleting $v_k$ would disconnect the nonempty sets
$A$ and $\{v_{k+1},\ldots,v_m\}$: no vertex of $A$ has an
off-path neighbour, and $v_0$ has no neighbour past $v_k$
by the definition of $k$.  These sets therefore could not be joined
through a vertex outside $K(P)$ either.  This contradicts
$2$-connectivity. Choose $j$ maximum and then $i$ minimum.

We check the availability of a representative for $v_iv_j$.
The path $Q_i$ traverses its initial block on $\{v_0,\ldots,v_k\}$
and then retains the original tail.  A tail representative containing \emph{both} $v_i$ and $v_j$ must
be $e_j$ or, when $j<m$, $e_{j+1}$: every other tail
representative has two prescribed key vertices distinct from
$v_i$ and $v_j$.
The latter cannot contain $v_i$: otherwise $v_iv_{j+1}$ would be a
crossing shadow edge, contrary to the maximal choice of $j$.
Any other representative from the re-ordered initial block has
two key vertices in $\{v_0,\ldots,v_k\}\setminus\{v_i\}$;
it cannot also contain both $v_i$ and $v_j$ in a triple.  Thus
only the first representative of $Q_i$ requires examination.
For $i\ge2$, the first representative joins $v_i$ to $v_{i-1}$;
if it also contains $v_j$, then $v_{i-1}v_j$ is a crossing shadow
edge, contradicting the minimal choice of $i$.
For $i=1$, its first representative is either an edge joining
$v_1$ to $v_0$, in which case it cannot contain $v_j$ since
$j>k$, or the exceptional edge $\{v_0,v_1,v_2\}$, which likewise
cannot contain $v_j$ as $j\ge4$.  Hence at most $e_j$ among
$E(Q_i)$ represents $v_iv_j$.  Because $v_i$ is low,
Lemma~\ref{G:lem:low-multiplicity} supplies a different
representative $f\in E(v_i,v_j)\setminus E(Q_i)$.
If $j=m$, this representative closes $Q_i$ to an $(m+1)$-cycle,
so $j<m$.

Apply the prefix-rotation ledger, Lemma~\ref{G:lem:prefix-rotation-ledger},
to $Q_i$ at $j$ with the outside representative $f$; this deletes
$e_j$ and inserts $f$. More precisely, reverse
the segment of $Q_i$ from its first vertex $v_i$ to $v_{j-1}$,
then traverse $f$ from $v_i$ to $v_j$ and retain the original tail.
The new path $Q_i^*$ uses exactly the following representative hyperedges
\[
 E(Q_i^*)=(E(Q_i)\setminus\{e_j\})\cup\{f\},
 \qquad f\notin E(Q_i).
\]
Thus no representative is repeated, and the resulting longest path has
endpoints $v_{j-1},v_m$. By the maximal endpoint-degree sum,
$\deg_H(v_{j-1})\le\deg_H(v_i)\le T$, so $v_{j-1}$ is low.
It shares the shadow neighbour $v_j$ with $v_i$ and hence is
adjacent to $v_i$ by \eqref{G:eq:low-common-neighbour}.
If $j\ge k+2$, then $v_0,v_{j-1}$ are distinct low vertices
with common neighbour $v_i$, contradicting the definition of $k$.
Consequently $j=k+1<m$.

Write $R$ for the rotated path.  Its endpoint $v_k$ is low, its
key set is $K(P)$, and its tail from $v_{k+1}$ onward is unchanged.
Here the rotation removes $e_{k+1}$ (the edge joining $v_k$ to
$v_{k+1}$ in $Q_i$) and inserts the outside witness $f$; in particular,
$R$ retains the original tail representatives $e_{k+2},\ldots,e_m$.
Also $N_H(v_k)\subseteq K(P)$.
Maximality of $j$ and the definition of $k$ imply that no vertex of $A$ has a neighbour
among $v_{k+2},\ldots,v_m$.  The $2$-connectivity of the shadow therefore supplies an edge
$v_kv_t$ with $t\ge k+2$. Otherwise deleting $v_{k+1}$ would
separate $\{v_0,\ldots,v_k\}$ from
$\{v_{k+2},\ldots,v_m\}$: the vertices in $A$ have no edges
to this tail, and the only remaining potential attachment from
the first set is $v_k$.  Here $v_k$ also has no off-path neighbour,
since it is a low endpoint of $R$; thus an off-path detour cannot
connect the two sets after deleting $v_{k+1}$. Choose $t$ maximum.

We next check the available representatives.
The path $R$ starts at $v_k$, follows the re-ordered initial block,
uses the representative $f\notin E(Q_i)$ to reach $v_{k+1}$,
and then follows the original tail.  Although $f$ need not lie
outside $E(P)$ (it could have been released in forming $Q_i$),
it occurs only once in $R$ and contains $v_i,v_{k+1}$.
Consequently $f$ cannot contain the pair $v_k,v_t$ for $t\ge k+2$.  A representative from the re-ordered block
other than the first has two vertices in $A$ or uses $f$, whose
two prescribed vertices are $v_i\in A$ and $v_{k+1}$;
none can contain both $v_k$ and $v_t$ with $t\ge k+2$.
The first representative has key vertices $v_k$ and a vertex of
$A$.  Were its third vertex $v_t$, that vertex of $A$ would be
adjacent to $v_t$, contrary to the preceding separation restriction.
In the unchanged tail, only $e_t$ and $e_{t+1}$ (if $t<m$)
can contain \emph{both} $v_k$ and $v_t$: any other tail
representative already has two prescribed key vertices different
from those two.  The latter cannot also contain $v_k$ by the
maximal choice of $t$.  Thus at most $e_t$ is a representative
of $R$ containing $v_k,v_t$; by low-degree multiplicity a second
one $h\in E(v_k,v_t)\setminus E(R)$ exists.  Notice that this
is a statement about \emph{the current path} $R$, not about the
original set $E(P)$; all representatives released by earlier
rotations have already been accounted for.  If $t=m$, the outside representative $h$ closes $R$, contrary to
the cycle barrier. Otherwise, apply Lemma~\ref{G:lem:prefix-rotation-ledger} to $R$
at $t$, reversing its initial segment up to $v_{t-1}$ and using $h$
to join $v_k$ to $v_t$.
The set of representative hyperedges is exactly
\[
 E(R^*)=(E(R)\setminus\{e_t\})\cup\{h\},
 \qquad h\notin E(R).
\]
This identity records the deletion of the old $v_{t-1}v_t$ incidence
and insertion of the new $v_kv_t$ incidence; hence no representative
is repeated. By maximality of the endpoint-degree sum,
$\deg_H(v_{t-1})\le\deg_H(v_k)\le T$, so $v_{t-1}$ is low.  The low vertices $v_k,v_{t-1}$ have common
neighbour $v_t$, so they are adjacent.
Now the low vertices $v_0,v_{t-1}$ have common neighbour $v_k$;
they too are adjacent, contradicting $t-1>k$.
\end{proof}

\begin{lemma}[Uniform cycle extension]\label{G:lem:uniform-cycle-extension}
Let $n\ge 9$, $r=\lfloor(n-1)/2\rfloor$, and $t=\binom r2$.
Let $H$ be a simple $3$-uniform hypergraph on $n$ vertices such that
\[
 d_H(u,v)=2\ \Longrightarrow\ \max\{\deg_H(u),\deg_H(v)\}\ge t+1.
\]
If $C$ is a Berge cycle of length $k<n$ and $x\notin K(C)$ satisfies
$N_H(x)\subseteq K(C)$ and $\deg_H(x)\ge t+1$, then $H$ has a Berge
cycle of length greater than $k$.
\end{lemma}
\begin{proof}
Suppose that no longer Berge cycle exists.  Among the cycles on
$K(C)$ choose one maximizing the number of representatives containing
$x$.  Let $s$ denote the number of vertices of $C$ incident with a
hyperedge through $x$ outside the representative set of $C$.
By Lemma~\ref{G:lem:cycle-count},
\begin{equation}\label{G:eq:unified-extension-bound}
 0\le s\le\lfloor k/2\rfloor,\qquad s\ne1,\qquad
 \deg_H(x)\le\binom s2+k-2s.
\end{equation}

First suppose $k\le n-2$ and $s\ge2$.  For $n\ge11$, the
convex function $f(s)=\binom s2+n-2-2s$ is at most $t$ on
$0\le s\le\lfloor(n-2)/2\rfloor$: at $s=0$ it is $n-2\le t$,
and at the other endpoint it is $\binom{r-1}2+1\le t$ for
$n=2r+1$, and $\binom r2=t$ for $n=2r+2$.
For $n\in\{9,10\}$ the same assertion on $s\ge2$ follows
by checking the endpoints of the intervals $2\le s\le3$ and
$2\le s\le4$, respectively.  Thus the bound
\eqref{G:eq:unified-extension-bound} contradicts $\deg_H(x)\ge t+1$.
If $k\le n-2$ and $s=0$, then $\deg_H(x)\le k\le n-2$; this
is impossible for $n\ge11$.  For $n=9,10$ we retain this case
and handle it below.

Now suppose $k=n-1$.  Set
$g(s)=\binom s2+n-1-2s$, where $0\le s\le r$.
For $2\le s\le r-1$, convexity and
\[
 g(2)=n-4\le t,\qquad
 g(r-1)=\begin{cases}
  \binom{r-1}2+2,&n=2r+1,\\
  \binom{r-1}2+3,&n=2r+2
 \end{cases}\le t
\]
rule out this range (recall $r\ge4$).  If $s=r$, then
$g(r)=t$ for odd $n$, and $g(r)=t+1$ for even $n$.
The latter equality case is treated below.
Finally $g(0)=n-1>t$ only for $n\in\{9,10,12\}$.

The numerical estimate has therefore isolated precisely two
possible configurations not already excluded by the degree threshold:
\begin{enumerate}[label=\textup{(\alph*)}]
\item $s=0$, with
\[
 (n,k)\in\{(9,7),(9,8),(10,7),(10,8),(10,9),(12,11)\};
\]
\item $n=2r+2$, $k=2r+1$, and $s=r$, with equality in
\eqref{G:eq:unified-extension-bound}.
\end{enumerate}
Indeed, in (a) the bound $\deg_H(x)\le k$ implies
$k\ge t+1$, and $k\le n-1$.  For odd $n=2r+1$, this is possible
only for $r=4$; for even $n=2r+2$, only for $r=4,5$.
The possible lengths are exactly those displayed above.

Consider (a). Every edge through $x$ is a representative of $C$:
otherwise an outside edge through $x$, whose two other vertices
belong to $K(C)$, would produce a nonempty index set in
Lemma~\ref{G:lem:cycle-count}.  Put
\[
 b=k-\deg_H(x),\qquad h=n-k-1.
\]
Here $b\ge0$, $k\ge\deg_H(x)\ge t+1\ge7$, and
\[
 b+h=n-1-\deg_H(x)\le n-t-2\le2.
\]
Hence $k>3b$ and
\[
 2+b+\binom{b+h}{2}\le2+2+\binom22=5\le t.
\]
Lemma~\ref{G:lem:few-exceptions}, with threshold $T=t$, produces
a $(k+1)$-cycle, a contradiction.

It remains to consider (b): $n=2r+2$, $k=n-1=2r+1$, and $s=r$.
Equality throughout \eqref{G:eq:unified-extension-bound} forces
exactly one cycle representative through $x$, and all
$\binom r2$ outside edges through $x$ on pairs of the $r$ indexed
vertices.  Relabel the cycle so that this representative is
$e_1=\{x,v_1,v_2\}$.  Lemma~\ref{G:lem:cycle-count} implies that
the index set is the unique independent $r$-subset of
$\{3,\ldots,2r+1\}$, namely $\{3,5,\ldots,2r+1\}$.
By Lemma~\ref{G:lem:parity-link}, with
$O=\{v_3,v_5,\ldots,v_{2r+1}\}$, we have
$N_H(v_q)\subseteq O$ for every even $q\in\{4,6,\ldots,2r\}$.
Consequently $\deg_H(v_4),\deg_H(v_6)\le\binom r2=t$.
These two vertices are nonadjacent, yet the consecutive cycle
representatives $e_4,e_5$ give a Berge $2$-path between them,
contradicting the distance-two hypothesis.
\end{proof}

\section{The case \texorpdfstring{$n\ge11$}{n >= 11}}\label{sec:large-order}
Put
\[
  r=\left\lfloor\frac{n-1}{2}\right\rfloor,
  \qquad t=\binom r2.
\]
Theorem~\ref{G:thm:low} handles longest paths with a low-degree
endpoint. We consider paths with both endpoint degrees at least $t+1$.

\begin{theorem}[Fan-type theorem]\label{L:thm:main}
Theorem~\ref{global:main} holds for $n\ge11$.
\end{theorem}

\begin{theorem}[High-degree endpoint case]\label{L:thm:high}
Let $H$ be a connected $3$-graph on $n\ge11$ vertices satisfying
Theorem~\ref{global:main}\textup{(a),(b)}.  If a longest Berge path
$P=v_0e_1v_1\cdots e_mv_m$ satisfies
$\deg_H(v_0),\deg_H(v_m)\ge t+1$, then $H$ contains a Hamiltonian Berge cycle.
\end{theorem}

\begin{proof}[Proof of Theorem~\ref{L:thm:main}]
Lemma~\ref{intro:edge-redundancy} shows that the hypotheses imply
that $H$ is $2$-edge-connected.
Choose a longest path maximizing the sum of its endpoint degrees.
Apply Theorem~\ref{L:thm:high} if both degrees exceed $t$, and
Theorem~\ref{G:thm:low} with $T=t$ otherwise.
\end{proof}

\subsection{Endpoint exchanges}
Suppose that $H$ and $P$ satisfy the hypotheses of
Theorem~\ref{L:thm:high}, and that $H$ has no Hamiltonian Berge cycle.
Lemma~\ref{G:lem:uniform-cycle-extension} supplies the length-$m$
extension property in Lemma~\ref{G:lem:high-structure}, with $T=t$.
Thus its boundary and crossing restrictions, and the charging
inequalities of Lemma~\ref{G:lem:endpoint-charging}, apply to $P$.
We need two further exchanges for the equality configurations.

\begin{lemma}[Five-slot saturation]\label{L:lem:five-slot}
Let $S\subseteq\{0,1,\ldots,m-1\}$ satisfy $|S|=5$ and $0\in S$.  Suppose
that all but at most one unordered pair from $\binom S2$ are occupied.  Assume
that an occupied pair $\{p,q\}$ is called \emph{usable at $p$} when
$|p-q|\ne1$.  Then among the indices
\[
 j\in S\setminus\{0,1\}
 \quad\text{for which }\{0,j\}\text{ is occupied},
\]
there exists a $j$ such that, after deleting the single pair $\{0,j\}$, the
vertex $j$ is still incident with an occupied usable pair.
\end{lemma}
\begin{proof}
At least two indices $j\in S\setminus\{0,1\}$ have
$\{0,j\}$ occupied. Call such a candidate bad if no occupied
nonconsecutive pair at $j$ remains after deleting $\{0,j\}$.
Among the three other indices, at most two are consecutive to $j$.
Thus a bad candidate forces the unique unoccupied pair to be
$\{j,q\}$ with $|j-q|\ne1$, and also forces $j-1,j+1\in S$.

If two candidates $j<\ell$ were bad, their forced unoccupied
pairs would coincide, so that pair would be $\{j,\ell\}$ and
$\ell-j\ge2$. Both intervals
$\{j-1,j,j+1\}$ and $\{\ell-1,\ell,\ell+1\}$ would lie in $S$.
Their union contains at least five positive indices, besides
$0\in S$, contrary to $|S|=5$. Hence at most one candidate is
bad, and a good candidate exists.
\end{proof}

\begin{lemma}[Saturated endpoint switch]\label{L:lem:saturated-switch}
Assume the hypotheses of Lemma~\ref{G:lem:high-structure}
with $T=t$ and $n\ge11$, let $m=2r+1$, and suppose
\[
 |I_1|=|I_2|=r,\qquad I_1=I_1',\qquad I_2=I_2',
\]
with a unique missing index $k\notin I_1\cup I_2$.  Assume $0<k<m-1$,
$k+1\in I_2$, and
\[
 e_1=\{v_0,v_1,v_{k+1}\}.
\]
If all $\binom r2$ pair-slots determined by $I_1$ are occupied by hyperedges
through $v_0$, then replacing $e_1$ by a suitable outside witness for
$\{v_0,v_1\}$ produces a longest Berge path $P^*$ with
\[
 |I_1(P^*)|\ge r+1,\qquad |I_2(P^*)|\ge r.
\]
Consequently its endpoint sets partition all $m$ indices.
\end{lemma}
\begin{proof}
Because $0\in I_1'=I_1$, choose
$a_0\in E(H)\setminus E(P)$ containing $\{v_0,v_1\}$ and exchange
it with $e_1$.  By Lemma~\ref{G:lem:first-edge-ledger}, this gives a
representative-safe $P^*$ with the same key order; the released $e_1$
witnesses $0,k\in I_1'(P^*)$.

Fix $i\in I_1\setminus\{0\}$.  For every $q\in I_1\setminus\{i\}$ the
occupied pair-slot $\{i,q\}$ gives
$h_{i,q}=\{v_0,v_{i+1},v_{q+1}\}$.  Choose $q$ so that
$|q-i|\ne1$ and $h_{i,q}\ne a_0$.  There are $r-1$ choices.  At most two are excluded by $|q-i|=1$, and at
most one further value can satisfy $h_{i,q}=a_0$, because for fixed $i$ the
equality $\{v_0,v_{i+1},v_{q+1}\}=a_0$ determines $v_{q+1}$, hence $q$,
uniquely.  Thus at most three choices are excluded; since $r\ge5$, one remains.
The two non-$v_0$ vertices are not consecutive on $P$, so $h_{i,q}$ is not a
path representative other than possibly $e_1$; but $e_1$ uses the missing
index $k$, whereas $i,q\in I_1$.
Thus $h_{i,q}\notin E(P^*)$ and still witnesses $i$.  The representative $h_{i,q}$ remains outside $E(P^*)$, and hence
$I_1(P^*)\supseteq I_1(P)\cup\{k\}$.

By Lemma~\ref{G:lem:high-structure}(ii), $0\notin I_2$.
Hence the outside edge $a_0$ cannot contain $v_m$: otherwise it would
witness $0\in I_2'$.  The first-representative exchange lemma thus
preserves every index in $I_2$, so $|I_2(P^*)|\ge r$.  Finally
Lemma~\ref{G:lem:endpoint-package}(1) gives
$|I_1(P^*)|+|I_2(P^*)|\le m=2r+1$, so equality holds.
\end{proof}

\subsection{Proof of the high-degree endpoint case}

\begin{proof}[Proof of Theorem~\ref{L:thm:high}]
Assume that $H$ has no Hamiltonian Berge cycle. For the path $P$ in
the statement, put $H'=H-E(P)$ and
\[
\begin{gathered}
 a=|I_1|,\quad a'=|I_1'|,\qquad
 b=|I_2|,\quad b'=|I_2'|,\\
 M=\{0,\ldots,m-1\}\setminus(I_1\cup I_2),\quad q=|M|.
\end{gathered}
\]
Lemma~\ref{G:lem:high-structure} supplies the disjointness,
boundary, and crossing restrictions.

\paragraph{Step 1: coarse bounds.}
The endpoint estimate gives
\[
 t+1\le\deg(v_0)\le\binom{a'}2+a+1
 \le\binom{a+1}2+1,
\]
so $a\ge r-1$; symmetrically $b\ge r-1$.

\paragraph{Step 2: refined bounds.}
Lemma~\ref{G:lem:endpoint-charging} gives
\begin{align}
 \deg(v_0)&\le\binom{a'}2+2(a-a')+q+1,
 \label{L:eq:refined-general-left}\\
 \deg(v_m)&\le\binom{b'}2+2(b-b')+q+1.
 \label{L:eq:refined-general-right}
\end{align}
If $a=r-1$, the coarse bound forces $a'=r-1$:
otherwise $\deg(v_0)\le\binom{r-2}2+r=t-r+3<t+1$.
Since $b\ge r-1$ and $m\le2r+1$, we have $q\le3$, whence
\[
 \deg(v_0)\le\binom{r-1}2+4=t-r+5\le t,
\]
a contradiction. Thus $a,b\ge r$, and therefore
\begin{equation}\label{L:eq:sum-near-full}
 2r\le a+b\le m\le2r+1,\qquad q\in\{0,1\}.
\end{equation}

\paragraph{Step 3: size-$r$ endpoint sets are prime.}
Suppose $a=r$. If $a'\le r-2$, the coarse bound gives
$\deg(v_0)\le\binom{r-2}2+r+1=t-r+4\le t$.
If $a'=r-1$, the refined bound and $q\le1$ give
$\deg(v_0)\le\binom{r-1}2+4=t-r+5\le t$.
Hence
\begin{equation}\label{L:eq:size-r-prime}
 |I_1|=r\quad\Longrightarrow\quad I_1=I_1',
\end{equation}
and the reversed statement holds for $I_2$.
If $q=1$, then $m=2r+1$ and $a=b=r$.
If $q=0$, at least one endpoint set has size $r$.

\paragraph{Case 1: $|I_1|+|I_2|=m$.}
At least one endpoint set has size $r$. Reverse the path if needed
so that $|I_1|=r$; Step~3 gives $I_1=I_1'$.
Lemma~\ref{G:lem:full-partition} then yields
$\deg_H(v_0)\le\binom r2=t$, a contradiction.
This argument excludes a full endpoint partition for \emph{every}
longest Berge path whose two endpoints have degree at least $t+1$;
we shall use this path-independent conclusion after a final rotation.

\paragraph{Case 2: $|I_1| + |I_2| = m-1$.}
Let $k$ be the unique index in
$\{0,1,\ldots,m-1\}\setminus(I_1\cup I_2)$. Since Step~2 gives
$|I_1|,|I_2|\ge r$ and $m\le n-1\le2r+1$, the equality
$|I_1|+|I_2|=m-1$ forces
\[
 m=2r+1=n-1,\qquad |I_1|=|I_2|=r.
\]
Step~3 gives $I_1=I_1'$ and $I_2=I_2'$.  Moreover
\eqref{L:eq:refined-general-left}--\eqref{L:eq:refined-general-right} with $q=1$
give
\[
 t+1\le\deg(v_0),\deg(v_m)\le t+2.
\]
There are five subcases. The boundary choices $k=0,m-1$ give
Subcases~2.1--2.2. Otherwise $0<k<m-1$ and both adjacent indices
lie in $I_1\cup I_2$. According as $k+1\in I_2$, or $k+1\in I_1$
with $k-1\in I_1$ or $I_2$, we obtain Subcases~2.3--2.5.

\paragraph{Subcase 2.1: $k=m-1$.}
Then $0\in I_1=I_1'$.  Let $a_0\in E(H')$ contain $\{v_0,v_1\}$.
By Lemma~\ref{G:lem:endpoint-package}(6), $N_H(v_0)\subseteq K(P)$, so write
$e_1=\{v_0,v_1,v_{j+1}\}$ with $1\le j\le m-1$.
If $j\in I_2$, Lemma~\ref{G:lem:first-edge-crossing} gives a
Berge cycle of length $m+1$.  If $j=k=m-1$, then
$e_1=\{v_0,v_1,v_m\}$ and
\[
 v_0\,a_0\,v_1e_2v_2\cdots e_mv_m\,e_1\,v_0
\]
is itself a Berge cycle of length $m+1$; the edge $a_0$ lies outside
$E(P)$ and all other representatives are the distinct path edges.
Thus $j\in I_1$.

There is also no missing-predecessor exception in the refined charging:
if a path edge $e_{i+1}$ containing $v_0$ has $i\in I_1'$ and
$i-1\notin I_2$, then $i-1$ cannot equal the missing index $m-1$ because
$i\le m-1$.  Hence $i-1\in I_1$.  We may therefore charge every edge
through $v_0$, including $e_1$ via $\{0,j\}$, injectively to a pair of
indices in $I_1$.  Consequently
\[
 \deg(v_0)\le\binom r2=t,
\]
a contradiction.

\paragraph{Subcase 2.2: $k=0$.}
This is the exact reverse-path analogue of Subcase~2.1.  Reversing a
Berge path changes only its order, not its representative set, so all
distinctness checks in that argument are preserved.  Here
$m-1\in I_2=I_2'$.  There is no missing-successor exception in the right
endpoint charging, and $e_m$ is either chargeable to a pair from $I_2$ or
directly produces, together with an outside witness for $m-1\in I_2'$, a
Berge cycle of length $m+1$.  Hence $\deg(v_m)\le t$, a contradiction.

\paragraph{Subcase 2.3: $0<k<m-1$ and $k+1\in I_2$.}
The missing-predecessor exception at the left endpoint is absent, because the
only possible exceptional index would be $i=k+1$, whereas
$k+1\notin I_1$.  Thus every edge through $v_0$ except possibly $e_1$ is
charged to an $I_1$-pair.  If $e_1$ is also chargeable, then
$\deg(v_0)\le\binom r2=t$, a contradiction.  If the third-vertex index
of $e_1$ belonged to $I_2$, Lemma~\ref{G:lem:first-edge-crossing}
would give a Berge cycle of length $m+1$.  Therefore an unchargeable $e_1$
must use the unique missing index $k$, and the only remaining possibility is
\begin{equation}\label{L:eq:23-exceptional-e1}
 e_1=\{v_0,v_1,v_{k+1}\},\qquad
 \deg(v_0)=t+1,
\end{equation}
and every one of the $\binom r2$ pair-slots determined by $I_1$ is occupied.

Lemma~\ref{L:lem:saturated-switch} applies to this equality configuration and
produces a longest Berge path $P^*$ with at least $r+1$ left endpoint indices
and at least $r$ right endpoint indices.  Since $m=2r+1$, the endpoint
disjointness lemma forces these sets to partition all $m$ indices.  Thus
$P^*$ lies in Case~1, already excluded.  This excludes Subcase~2.3.

\paragraph{Subcase 2.4: $0<k<m-1$, $k+1\in I_1$, and $k-1\in I_1$.}
Reverse the path.  Under the reversal correspondence
\[
 i\longmapsto m-1-i,
 \qquad I_1(P)\longleftrightarrow I_2(\overleftarrow P),
 \qquad I_2(P)\longleftrightarrow I_1(\overleftarrow P),
\]
the unique missing index $k$ becomes $m-1-k$, and the assumption
$k-1\in I_1(P)$ becomes exactly the successor condition required in
Subcase~2.3 for the reversed path.  Hence that subcase applies.

\paragraph{Subcase 2.5: $0<k<m-1$, $k+1\in I_1$, and $k-1\in I_2$.}
We first deal with the endpoint representatives. Write $S=I_1=I_1'$, so
$|S|=r$ and $|\binom S2|=t$. The crossing restriction and the unique
missing index imply that every path representative through $v_0$, except
$e_1$ and possibly $e_{k+2}$, is charged to a pair-slot in $\binom S2$:
for $e_{i+1}$ with $i\in S$ and $i\ge1$, the predecessor $i-1$ lies in $S$ unless
$i=k+1$, which is precisely the exceptional slot $\mu_k$. Outside edges
through $v_0$ are charged to their pair of indices in $S$. If $e_1$ is not
pair-chargeable, it occupies the separate slot $\star_0$, while $e_{k+2}$
is the only possible preimage of $\mu_k$. Consequently
\[
 \deg(v_0)\le |\mathcal P_{\rm occ}|+2,
 \qquad \mathcal P_{\rm occ}\subseteq\binom S2,
\]
where $\mathcal P_{\rm occ}$ is the set of occupied pair-slots. Since
$\deg(v_0)\ge t+1$, at least $t-1$ of the $t$ pair-slots are occupied;
there is at most one unoccupied pair-slot. As in Case~1, a third-vertex
index from $I_2$ would give a Berge cycle by
Lemma~\ref{G:lem:first-edge-crossing}. Therefore, if $e_1$ is not
pair-chargeable, its third-vertex index must be the unique missing index
$k$, that is, $e_1=\{v_0,v_1,v_{k+1}\}$.
Choose $a_0\in E(H')$ with $\{v_0,v_1\}\subseteq a_0$ and replace $e_1$
by $a_0$.  As in Subcase~2.3, the old $e_1$ makes $k$ an outside index for
the new path.  We choose the replacement $a_0$ so that no old member
of $I_1$ is lost.  For $r\ge6$, the previous counting suffices: for each
fixed $i\in I_1\setminus\{0\}$, among the $r-1$ pair-slots containing $i$,
at most one is globally unoccupied, at most two correspond to consecutive
non-$v_0$ vertices, and at most one is the newly used edge $a_0$; hence an
occupied nonconsecutive slot remains.

When $r=5$, put $S:=I_1$.  At most one pair of $\binom S2$ is
unoccupied.  Lemma~\ref{L:lem:five-slot} gives an index
$j\in S\setminus\{0,1\}$ such that $\{0,j\}$ is occupied and, after this
pair is removed, $j$ still has an occupied nonconsecutive pair-slot.  Let
\[
 a_0=a_j:=\{v_0,v_1,v_{j+1}\}
\]
be the outside edge occupying $\{0,j\}$.  This edge is indeed outside
$E(P)$: a path representative charged to a pair-slot has the consecutive form
$\{i-1,i\}$, while $j\ne1$.  No old left endpoint index is lost.  The index
$0$ is witnessed by the old edge $e_1$, which is outside the new path.  The
good index $j$ retains, by Lemma~\ref{L:lem:five-slot}, an occupied usable
pair-slot after $\{0,j\}$ is removed.  Finally, if
$i\in I_1\setminus\{0,j\}$, the newly used pair $\{0,j\}$ is not incident
with $i$.  Among the four pair-slots of $S$ incident with $i$, at most one is
globally unoccupied and at most two correspond to consecutive non-$v_0$
vertices.  Hence at least one occupied nonconsecutive pair-slot remains.
Thus every old index $i\in I_1$ is preserved: $0$ by the old $e_1$, $j$ by
the five-slot lemma, and every other index by the preceding count.  For
$i\ne0$, the hyperedge occupying the chosen nonconsecutive slot cannot be
$e_1$, because $e_1$ uses the unique missing index $k\notin I_1$, whereas
both indices of the chosen pair lie in $I_1$.  Hence this hyperedge is outside
the new path and still witnesses $i$.  As in Subcase~2.3,
removing $a_0$ from the outside-edge set cannot destroy an $I_2'$-witness
unless $a_0$ also contains $v_m$, in which case $a_0$ closes $P$ and
already gives a contradiction.  Hence the new path has at least $r+1$
left indices and at least $r$ right indices, so it lies in Case~1, a
contradiction.  Therefore $e_1$ is chargeable to an $I_1$-pair.  The
reverse argument shows that $e_m$ is chargeable to an $I_2$-pair.

We can now invoke the exact charging map rather than an
implicit equality argument.  At the left endpoint, the hypotheses
of Lemma~\ref{G:lem:one-gap-rigidity} hold with $S=I_1$;
at the right endpoint they hold after path reversal.
Since both endpoint degrees exceed $t=\binom r2$, the two
finite injections are bijective.  In particular,
\begin{equation}\label{L:eq:25-left-exact}
\deg(v_0)=t+1,\qquad
e_{k+2}=\{v_0,v_{k+1},v_{k+2}\},
\end{equation}
\begin{equation}\label{L:eq:25-right-exact}
\deg(v_m)=t+1,\qquad
e_k=\{v_m,v_{k-1},v_k\}.
\end{equation}
Lemma~\ref{G:lem:one-gap-rigidity} also gives
\[
 \deg_H(v_0,v_{k+1})=\deg_H(v_m,v_k)=1.
\]
Condition~\textup{(b)} therefore yields
$\deg_H(v_{k+1}),\deg_H(v_k)\ge t+1$.

Choose $b\in E(H')$ containing $\{v_m,v_{k-1}\}$; such an edge exists
because $k-1\in I_2'=I_2$.  Delete $e_k$ from $P$ and insert $b$, obtaining
the Hamilton Berge path
\begin{equation}\label{L:eq:25-rotated-path}
 Q=
 v_k e_{k+1}v_{k+1}e_{k+2}\cdots e_mv_m\,
 b\,v_{k-1}e_{k-1}\cdots e_1v_0.
\end{equation}
Its representative set is
$(E(P)\setminus\{e_k\})\cup\{b\}$, so all representatives are distinct.
The endpoints $v_k,v_0$ both have degree at least $t+1$.

Apply only the numerical estimates of Steps~1--3 to $Q$, redefining
the endpoint index sets relative to the displayed order of $Q$.
These steps use the endpoint-package and charging lemmas, not the
case analysis that follows them.  The path $Q$ is longest
and both of its endpoints have degree at least $t+1$, so the same structural
corollary and counts apply.  Case~1 has already been ruled out for every
longest path with high-degree endpoints.  Hence $Q$ must have
exactly one missing endpoint index, and Step~3 then says that each of its two
size-$r$ endpoint sets equals its outside-witness subset. This contradicts primeness at the right endpoint $v_0$.
Indeed, in the order displayed in
\eqref{L:eq:25-rotated-path}, put $w_0=v_k,w_1=v_{k+1},w_2=v_{k+2}$.
The complete representative ledger for the rotation is
\[
 E(Q)=\bigl(E(P)\setminus\{e_k\}\bigr)\cup\{b\},
 \qquad b\notin E(P),
\]
with $e_{k+2}$ retained in its original representative set;
no edge is counted twice. The second path representative is
\[
 e_{k+2}=\{v_0,w_1,w_2\},
\]
so index $1$ belongs to the right endpoint set of $Q$.  But
$\deg(v_0,w_1)=\deg(v_0,v_{k+1})=1$, and its unique hyperedge
$e_{k+2}$ is already a representative of $Q$.  Thus index $1$ has no
outside witness, contradicting the conclusion of Step~3.  This excludes Subcase~2.5.

The five cases exhaust the possibilities.
\end{proof}

\section{The boundary cases \texorpdfstring{$n=9$ and $n=10$}{n = 9 and n = 10}}\label{sec:boundary}
For $n\in\{9,10\}$ we have $r=4$ and $t=6$.  The two local hypotheses become
\begin{align}
 d_H(u,v)=2 &\Longrightarrow \max\{\deg_H(u),\deg_H(v)\}\ge7, \tag{F}\label{B:cond:F}\\
 \deg_H(u,v)=1 &\Longrightarrow \min\{\deg_H(u),\deg_H(v)\}\ge7. \tag{C}\label{B:cond:C}
\end{align}
We use the lemmas of Section~\ref{sec:common-tools}, independently
of the proof for $n\ge11$.

\begin{theorem}[Orders nine and ten]\label{B:thm:main}
Let $H$ be a $2$-edge-connected simple $3$-uniform hypergraph on $n\in\{9,10\}$ vertices whose $2$-shadow is $2$-connected. Put
\[
 r=\left\lfloor\frac{n-1}{2}\right\rfloor=4,
 \qquad t=\binom42=6.
\]
Assume that every pair at Berge distance two has at least one endpoint of degree at least $7$, and that every pair of codegree one has both endpoints of degree at least $7$. Then $H$ contains a Hamiltonian Berge cycle.
\end{theorem}

\subsection{Endpoint counts and exchanges}
Throughout this subsection, assume that $H$ is a connected $3$-graph
of order $n\in\{9,10\}$ satisfying \eqref{B:cond:F}--\eqref{B:cond:C},
that $H$ has no Hamiltonian Berge cycle, and that
$P=v_0e_1v_1\cdots e_mv_m$ is a longest Berge path with both endpoint
degrees at least seven. Lemma~\ref{G:lem:uniform-cycle-extension}
supplies the extension hypothesis of Lemma~\ref{G:lem:high-structure}
with $T=6$. In particular its boundary and crossing restrictions apply.
Write
\[
 a=|I_1|,\quad a'=|I_1'|,\quad b=|I_2|,\quad b'=|I_2'|,
 \qquad M=\{0,\ldots,m-1\}\setminus(I_1\cup I_2),\quad q=|M|.
\]
The common endpoint and charging lemmas give
\begin{equation}\label{B:eq:coarse}
 \deg_H(v_0)\le\binom{a'}2+a+1,\qquad
 \deg_H(v_m)\le\binom{b'}2+b+1,
\end{equation}
and
\begin{align}
 \deg_H(v_0)&\le\binom{a'}2+2(a-a')+q+1,
 \label{B:eq:refined-left}\\
 \deg_H(v_m)&\le\binom{b'}2+2(b-b')+q+1.
 \label{B:eq:refined-right}
\end{align}
We use the injection $\Phi_0$ from
Lemma~\ref{G:lem:endpoint-charging}; equality makes it bijective.
If the endpoint sets partition all indices and a side of size four
equals its prime subset, Lemma~\ref{G:lem:full-partition} bounds its
degree by six, a contradiction. The lemmas below address the
remaining order-ten configurations.

\begin{lemma}[Four-index safe replacement]\label{B:lem:safe-replacement}
Assume the hypotheses of \cref{G:lem:high-structure} with $T=6$ and $n=10$, and let
$P=v_0e_1v_1\cdots e_mv_m$ be a longest Berge path in the central
high-endpoint configuration $m=9$ with a unique missing index
$0<k<m-1$.  Suppose
\[
 S:=I_1=I_1',\qquad |S|=4,\qquad 0\in S,
\]
where $k\notin I_1\cup I_2$ is the unique missing endpoint index. Assume that at least five of the six pair-slots in $\binom S2$ are occupied under the pair-slot part of the injection $\Phi_0$ in \cref{G:lem:endpoint-charging}, and that
\[
 e_1=\{v_0,v_1,v_{k+1}\}.
\]
Then there is an outside edge $a_0\in E(H)\setminus E(P)$ containing $v_0,v_1$ such that, after replacing $e_1$ by $a_0$, every old member of $I_1$ is still a left endpoint index and $k$ becomes an outside-witness index. Moreover every old right endpoint index is preserved. Consequently, if $P^*$ denotes the new path, then
\[
 |I_1(P^*)|\ge5,
 \qquad
 I_2(P)\subseteq I_2(P^*).
\]
\end{lemma}

\begin{proof}
Because $0\in I_1'$, at least one outside edge contains $v_0,v_1$. Any such edge has the form
\[
 a_j=\{v_0,v_1,v_{j+1}\}
\]
for some $j\in S\setminus\{0\}$; call it an outside candidate. Under the explicit charging of \cref{G:lem:endpoint-charging}, an occupied pair-slot $\{0,j\}$ with $j\ne1$ must be occupied by this outside edge: a path representative can occupy a pair-slot only when that pair is consecutive, and the only path representative that can occupy $\{0,1\}$ is $e_2=\{v_0,v_1,v_2\}$. Thus the only exceptional $0$-slot is $\{0,1\}$.

Call an outside candidate $a_j$ \emph{bad} if, after $a_j$ is moved into the path in place of $e_1$, the old index $j$ has neither an outside witness nor a special path-edge witness. Only indices $0$ and $j$ can lose an outside witness, by the precise outside-set identity in \cref{G:lem:first-edge-ledger}.

A bad candidate has the following form. Besides $0,j$, the set $S$ has exactly two indices. Any occupied pair-slot $\{j,q\}$ with $|q-j|\ne1$ is necessarily occupied by an outside edge through $v_0$, and that edge continues to witness $j$. If the predecessor slot $\{j-1,j\}$ is occupied, then either it is occupied by an outside edge, which again witnesses $j$, or it is occupied by the path representative
\[
 e_{j+1}=\{v_0,v_j,v_{j+1}\},
\]
which is precisely the special path-edge witness for $j$. Hence a bad $j$ is possible only when $j+1\in S$, the successor slot $\{j,j+1\}$ is occupied by the successor path representative
\[
 e_{j+2}=\{v_0,v_{j+1},v_{j+2}\},
\]
which does not witness $j$, and the remaining non-$0$ pair-slot incident with $j$ is the unique unoccupied pair-slot. Thus a bad candidate determines the unique unoccupied pair.

At most one outside candidate is bad. Suppose $j<\ell$ are both bad. Then $j+1,\ell+1\in S$. If $\ell\ge j+2$, the five indices
\[
 0,j,j+1,\ell,\ell+1
\]
are distinct, contradicting $|S|=4$. Hence $\ell=j+1$. Since $j$ is bad, the pair $\{j,j+1\}$ is occupied. For $\ell=j+1$ this is the predecessor slot $\{\ell-1,\ell\}$; if it is occupied by an outside edge, it witnesses $\ell$, while if it is occupied by a path representative, that representative is $e_{\ell+1}=e_{j+2}$ and is the special witness for $\ell$. Either way $\ell$ is not bad, contradiction.

If there are at least two outside candidates, one is therefore good. Suppose there is only one. Then $1\in S$; otherwise all three pairs $\{0,q\}$ with $q\in S\setminus\{0\}$ are nonconsecutive, and since at most one pair-slot is unoccupied, at least two would be outside candidates. With $1\in S$, the unique candidate cannot be $\{0,1\}$, for then the two remaining $0$-slots are nonconsecutive and at most one is unoccupied, so the other would give a second outside candidate. Hence among the two $0$-slots other than the unique outside candidate, one is the consecutive slot $\{0,1\}$ and the other is the unique unoccupied pair. Thus the unique unoccupied pair contains $0$, whereas badness of the sole candidate would force the unique unoccupied pair to be a non-$0$ pair incident with that candidate. Therefore the sole candidate is good.

Choose a good outside candidate $a_0=a_j$ and replace $e_1$ by $a_0$.
The key-vertex order is unchanged and the new representative set is
$(E(P)\setminus\{e_1\})\cup\{a_0\}$.  The old edge $e_1$, now outside the
new path, contains both $\{v_0,v_1\}$ and $\{v_0,v_{k+1}\}$, so it witnesses
the old index $0$ and the missing index $k$. By goodness, $j$ retains a
representative-safe witness, and no other old member of $S$ can lose one.
Hence every old member of $I_1$ is preserved and $k$ is added.

Since $j\in I_1$ but $m-1\notin I_1$, we have $j+1\ne m$,
so $v_m\notin a_0$.  The exchange identity in
Lemma~\ref{G:lem:first-edge-ledger} therefore gives
$I_2(P)\subseteq I_2(P^*)$ and verifies distinctness of all
representatives of $P^*$.
\end{proof}

\begin{lemma}[The $(4,3)$ critical endpoint]\label{B:lem:critical}
Assume the hypotheses and notation of \cref{G:lem:high-structure} and
\cref{G:lem:endpoint-charging}, with $n=10$, $r=4$, $T=t=6$, $q=1$, and
\[
 |I_1|=4,
 \qquad |I_1'|=3,
 \qquad \deg_H(v_0)\ge7.
\]
Let $p$ be the unique element of $I_1\setminus I_1'$ and let $k$ be the unique missing endpoint index. Then equality holds in \eqref{B:eq:refined-left}; hence $\Phi_0$ is bijective. In particular
\begin{equation}\label{B:eq:critical-rel}
 p+1\in I_1',\qquad k+1\in I_1',
\end{equation}
and every pair-slot in $\binom{I_1'}2$ is occupied.
\end{lemma}

\begin{proof}
Substituting $a=4$, $a'=3$, $q=1$ into \eqref{B:eq:refined-left} gives
\[
 7\le \deg_H(v_0)\le \binom32+2(1)+1+1=7.
\]
Thus equality holds and, by the explicit injection of \cref{G:lem:endpoint-charging}, $\Phi_0$ is a bijection onto
\[
 \binom{I_1'}2\,\dot\cup\,\{D_p\}\,\dot\cup\,\{R_p\}\,\dot\cup\,\{M_k\}\,\dot\cup\,\{\star_0\}.
\]
Therefore every pair-slot is occupied. Occupancy of the predecessor slot $R_p$ means that some prime index $i\in I_1'$ with $i-1=p$ has its path representative through $v_0$; necessarily $i=p+1$, proving $p+1\in I_1'$. Occupancy of the missing-predecessor slot $M_k$ similarly means that a prime index $i\in I_1'$ satisfies $i-1=k$, so $i=k+1$ and $k+1\in I_1'$. This proves \eqref{B:eq:critical-rel}.
\end{proof}

\begin{lemma}[Normalization of a $(4,3)$ endpoint]\label{B:lem:normalization}
Assume the hypotheses and notation of \cref{G:lem:high-structure} with
$T=6$ and $n=10$. Suppose
\[
 m=9,\qquad |I_1|=|I_2|=4,
\]
there is a unique missing endpoint index $k\notin I_1\cup I_2$, and
\[
 0<k<m-1,
 \qquad k+1\in I_1,
 \qquad k-1\in I_2.
\]
If the left endpoint has critical type $(4,3)$, one can exchange the
first representative to obtain one of the following outcomes: a Berge cycle
of length $m+1$; a longest Berge path $P^*$ whose endpoint index sets
partition $\{0,\ldots,m-1\}$; or a longest Berge path $P^*$ with the same
key-vertex order and endpoint degrees for which
$ I_1(P^*)=I_1'(P^*)$ and the unique missing index remains $k$.
In either path outcome the exchange preserves
\[
 I_2(P)\subseteq I_2(P^*),\qquad I_2'(P)\subseteq I_2'(P^*).
\]
\end{lemma}

\begin{proof}
Let $p$ be the unique member of $I_1\setminus I_1'$. Since $0\notin I_2$ and $k$ is the unique missing index, $0\in I_1$. Index $0$ cannot be path-only, because that alternative would require $e_1=\{v_0,v_0,v_1\}$. Hence $0\in I_1'$. By \cref{B:lem:critical},
\[
 I_1'=\{0,p+1,k+1\}.
\]
Bijectivity of $\Phi_0$ also identifies the path representatives occupying the predecessor and missing-predecessor slots:
\[
 e_{p+2}=\{v_0,v_{p+1},v_{p+2}\},
 \qquad
 e_{k+2}=\{v_0,v_{k+1},v_{k+2}\}.
\]
The central relations force $p\notin\{k-1,k,k+1\}$: indeed,
$p\in I_1$, $k-1\in I_2$, $k$ is missing, and
$k+1\in I_1'$.  Thus $|p-k|\ge2$.  In particular, the
three prime indices $0,p+1,k+1$ are pairwise nonconsecutive.
The bijective charging map now determines the \emph{entire}
left-endpoint star.  Set
\[
 A_p=\{v_0,v_1,v_{p+2}\},\qquad
 A_k=\{v_0,v_1,v_{k+2}\},\qquad
 h=\{v_0,v_{p+2},v_{k+2}\}.
\]
Since all three pair-slots of $\binom{I_1'}2$ are occupied and
none is consecutive, their occupants are exactly the three outside
triples $A_p,A_k,h$.  The four remaining slots
$\star_0,D_p,R_p,M_k$ are occupied by
$e_1,e_{p+1},e_{p+2},e_{k+2}$, respectively.  Consequently
\begin{equation}\label{B:eq:critical-star}
 E_H(v_0)=
 \{A_p,A_k,h,e_1,e_{p+1},e_{p+2},e_{k+2}\},
 \qquad A_p,A_k,h\notin E(P).
\end{equation}
This is an equality of \emph{actual triples}, not merely of
charging slots.  The path representatives and outside triples are
pairwise distinct, as required by the injection.

Choose the explicit outside candidate $a_0=A_p$ and replace
$e_1$ by $a_0$.  Since $p+1\in I_1'$ but $m-1\notin I_1$,
we have $p+2<m$, so $v_m\notin a_0$.  By
\cref{G:lem:first-edge-ledger}, only the old outside-witness indices
$0$ and $p+1$ might be affected.  The released $e_1$ witnesses
$0$, while $h$, which contains no $v_1$ and is distinct from $a_0$,
remains outside the new path and witnesses both $p+1$ and $k+1$.
Thus all old prime indices survive.  The sole old path-only index
$p$ retains its path representative $e_{p+1}$.  In particular,
\[
 I_1'(P)\subseteq I_1'(P^*),\qquad
 I_1(P)\subseteq I_1(P^*).
\]
Write
\[
 e_1=\{v_0,v_1,v_{\ell+1}\}.
\]
The index $\ell$ cannot lie in $I_2$, by \cref{G:lem:first-edge-crossing}; this crossing principle does not require the full-partition hypothesis. We also claim $\ell\notin I_1'$. If $\ell\in I_1'$, then the pair-slot $\{0,\ell\}\in\binom{I_1'}2$ is occupied because $\Phi_0$ is bijective. Any edge occupying this pair-slot is necessarily the triple $\{v_0,v_1,v_{\ell+1}\}=e_1$. But in the bijection $e_1$ is the unique preimage of the disjoint slot $\star_0$, so it cannot also occupy $\{0,\ell\}$; a second edge with the same triple is forbidden by simplicity. Thus $\ell\notin I_1'$.

The indices form the disjoint union
\[
 I_1'\,\dot\cup\,\{p\}\,\dot\cup\,I_2\,\dot\cup\,\{k\},
\]
so $\ell\in\{p,k\}$. If $p=1$, the alternative $\ell=p$ is impossible:
the path-only index $1$ gives $e_2=\{v_0,v_1,v_2\}$, while
$\ell=p=1$ would give $e_1=\{v_0,v_1,v_2\}$, contradicting
distinctness of representatives. Thus $p=1$ necessarily falls
into the $\ell=k$ outcome. If $\ell=k$, the released $e_1$ witnesses the previously missing
index $k$. The new left endpoint set contains at least five indices;
together with the four preserved right indices, the disjointness of
\cref{G:lem:endpoint-package} forces a full partition of the nine
indices. If $\ell=p$, the released $e_1$ makes the unique path-only
index $p$ prime. The left endpoint set is then unchanged and satisfies
$I_1(P^*)=I_1'(P^*)$.

Finally $v_m\notin a_0$, so the exchange removes no outside
witness at the opposite endpoint.  Since $e_1$ does not contain
$v_m$ in the central configuration, no right path-only witness is
lost either. Thus
\[
 I_2'(P)\subseteq I_2'(P^*),\qquad
 I_2(P)\subseteq I_2(P^*).
\]
In the prime,
non-full-partition outcome the left set remains the same four indices; since
the right set cannot increase without producing a full partition, the unique
missing index remains $k$.  The new representative set is
$(E(P)\setminus\{e_1\})\cup\{a_0\}$, so all representatives are pairwise
distinct.
\end{proof}

\begin{lemma}[Two-sided normalization preserves established primeness]
\label{B:lem:two-sided-normalization}
Assume the central $n=10$ configuration of
Lemma~\ref{B:lem:normalization}. Assume, in addition, that
\begin{enumerate}[label=\textup{(N\arabic*)}]
\item no longest Berge path whose two endpoints have degree at least $7$
has endpoint index sets that partition all nine indices; and
\item every such path has both endpoint index sets of size at least four.
\end{enumerate}
Then one may normalize the left endpoint, if critical, and then the
right endpoint, if critical, without destroying a prime side already
obtained. Unless a forbidden Berge cycle or full partition occurs, the
same missing index $k$ remains, both endpoint sets have size four, and
$I_1=I_1'$, $I_2=I_2'$ after at most two exchanges.
\end{lemma}
\begin{proof}
Each application of Lemma~\ref{B:lem:normalization}, including its
reversed form, preserves the opposite endpoint's whole index set and its
prime subset. For two successive exchanges, denoting the initial path by
$P$, the left-normalized path by $P^L$, and the subsequent
right-normalized path by $P^{LR}$, the ledger specifically gives
\[
 I_1'(P^L)\subseteq I_1'(P^{LR}),\qquad
 I_2'(P)\subseteq I_2'(P^L),
\]
with index labels read in the original key order. The first inclusion
records the prime indices protected by the second (reversed) exchange;
the second records those protected by the first exchange. The two hypotheses (N1)--(N2) continue to apply because each
exchange keeps the same key vertices and endpoint degrees. If the left
side is already prime, no first exchange is needed; otherwise apply
Lemma~\ref{B:lem:normalization}. Unless this produces one of the two
contradictory outcomes, the resulting left set equals its prime subset
and has size four. The right set contains its four original indices.
Since the sets are disjoint and a full partition is excluded by (N1),
neither set can grow beyond four; thus both have size exactly four, the
right set is unchanged, and the single missing index remains $k$.

If the right side is not prime, apply the reversed normalization.
It preserves all four \emph{prime} left indices, so after this step the
left prime subset still has cardinality four. By (N1)--(N2) and
endpoint-set disjointness, neither endpoint set can have more than four
indices; hence the left set still equals its preserved prime subset, and
the right side has just become prime.  The orientation reversal changes the index
labels by $i\mapsto m-1-i$, but returning to the original order
restores the same key-vertex order and missing index $k$.
Every exchange changes exactly one representative at a time, as
recorded in Lemma~\ref{G:lem:first-edge-ledger}, so representative
distinctness is maintained at both stages.
\end{proof}

\begin{theorem}[High-degree endpoint theorem for orders nine and ten]\label{B:thm:high}
Let $H$ be a connected $3$-graph on $n\in\{9,10\}$ satisfying \eqref{B:cond:F} and \eqref{B:cond:C}. If a longest Berge path
\[
 P=v_0e_1v_1\cdots e_mv_m
\]
satisfies $\deg_H(v_0),\deg_H(v_m)\ge t+1$, then $H$ contains a Hamiltonian Berge cycle.
\end{theorem}

\begin{proof}
Assume, for contradiction, that $H$ contains no Hamiltonian Berge cycle.
Then \cref{G:lem:high-structure,G:lem:endpoint-charging} and the notation above apply to $P$.

\emph{Order $n=9$.}
Here $r=4$, $t=6$, and $m\le8$. From \eqref{B:eq:coarse}, $a,b\ge3$. Suppose $a=3$. The coarse bound and $\deg_H(v_0)\ge7$ force $a'=3$. Then \eqref{B:eq:refined-left} gives
\[
 7\le3+q+1,
\]
so $q\ge3$. But $b\ge3$ and $m\le8$ give
\[
 q=m-a-b\le8-3-3=2,
\]
contradiction. Thus $a,b\ge4$, and hence
\[
 8\le a+b\le m\le8.
\]
Therefore $m=8$, $a=b=4$, and $q=0$. If $a'\le2$, \eqref{B:eq:coarse} gives $\deg_H(v_0)\le6$; if $a'=3$, \eqref{B:eq:refined-left} with $q=0$ gives
\[
 \deg_H(v_0)\le\binom32+2+1=6.
\]
Thus $a'=4$, and similarly $b'=4$. The endpoint sets partition all eight indices and both size-four sides are prime. \Cref{G:lem:full-partition} gives a contradiction.

\emph{Order $n=10$.}
Again $r=4$, $t=6$, now $m\le9$. First rule out $a=3$. The coarse bound forces $a'=3$ and $\deg_H(v_0)=7$. Then \eqref{B:eq:refined-left} gives $q\ge3$. Since $b\ge3$ and $m\le9$, equality holds throughout:
\[
 m=9,
 \qquad a=b=3,
 \qquad a'=b'=3,
 \qquad q=3.
\]
Let $M$ be the three missing indices. At the left endpoint, the target of $\Phi_0$ has size
\[
 \binom32+3+1=7,
\]
so $\Phi_0$ is bijective. Thus every missing-predecessor slot $M_h$ is occupied, forcing
\[
 h+1\in I_1'\qquad(h\in M).
\]
By the symmetric bijection at the right endpoint,
\[
 h-1\in I_2'\qquad(h\in M).
\]
Hence $M\subseteq\{1,\ldots,7\}$. Two missing indices cannot differ by one, because the earlier one would have its successor in $I_1'$; they cannot differ by two, because the intermediate index would belong to both $I_1$ and $I_2$. Thus three elements of $\{1,\ldots,7\}$ with pairwise gaps at least three must be $\{1,4,7\}$. But $7\in M$ forces $8\in I_1'$, contradicting $m-1=8\notin I_1$. Therefore $a\ge4$, and symmetrically $b\ge4$.

If $q=0$, then $a+b=m\le9$, so at least one endpoint set has size four. For a size-four left set, $a'\le2$ gives the coarse bound $\deg_H(v_0)\le6$, while $a'=3$ gives
\[
 \deg_H(v_0)\le\binom32+2+1=6.
\]
Thus a size-four high-degree endpoint is prime; the right statement is symmetric. \Cref{G:lem:full-partition} gives a contradiction.

It remains to consider $q=1$. Then $m=9$ and $a=b=4$. The refined bound shows that a size-four high-degree endpoint set is either prime or has the unique critical type $(4,3)$ of \cref{B:lem:critical}. Let $k$ be the unique missing index.

First eliminate the noncentral locations. If $k=m-1$, then $0\in I_1$. The left endpoint cannot be critical, because \cref{B:lem:critical} would require $k+1=m\in I_1'$, outside the index range. Hence $I_1=I_1'$. There is no missing-predecessor exception at the left endpoint. By the charging argument in the proof of \cref{G:lem:full-partition}, every edge through $v_0$ is chargeable to an $I_1$-pair unless $e_1=\{v_0,v_1,v_m\}$, in which case an outside witness for $0\in I_1'$ closes a Berge cycle of length $m+1$. Hence $\deg_H(v_0)\le6$, contradiction. The case $k=0$ is symmetric.

Assume $0<k<m-1$. Both $k-1$ and $k+1$ lie in $I_1\cup I_2$. If $k+1\in I_2$, the left endpoint cannot be critical, because \cref{B:lem:critical} would require $k+1\in I_1'\subseteq I_1$. Thus $I_1=I_1'$. Since the only possible missing-predecessor contribution would occur at a prime index $i=k+1$, and here $k+1\notin I_1$, every edge through $v_0$ except possibly $e_1$ is charged to an $I_1$-pair. If $e_1$ is pair-chargeable, then $\deg_H(v_0)\le6$. Otherwise its third-vertex index cannot lie in $I_2$ by \cref{G:lem:first-edge-crossing}; hence it is $k$, so
\[
 e_1=\{v_0,v_1,v_{k+1}\}.
\]
Since $\deg_H(v_0)\ge7$, at least five of the six $I_1$-pair slots are occupied. Let $P^*$ be the longest path supplied by \cref{B:lem:safe-replacement}. That lemma gives $|I_1(P^*)|\ge5$ and preserves all four old right indices, so $|I_2(P^*)|\ge4$. The endpoint-package crossing restriction applied to $P^*$ gives $I_1(P^*)\cap I_2(P^*)=\varnothing$, while there are only $m=9$ endpoint indices. Hence equality holds throughout:
\[
 |I_1(P^*)|=5,\qquad |I_2(P^*)|=4,
\]
and the two endpoint sets partition all nine indices. Since the right endpoint still has degree at least seven, the $q=0$ refined endpoint count used above forces its size-four endpoint set to be prime. Thus \cref{G:lem:full-partition} applies, a contradiction. The case in which both $k+1$ and $k-1$ lie in $I_1$ is the reverse analogue.

Thus the only remaining location is the central subcase
\begin{equation}\label{B:eq:central}
 0<k<m-1,
 \qquad k+1\in I_1,
 \qquad k-1\in I_2.
\end{equation}
The boundary exclusions $0\notin I_2$ and $m-1\notin I_1$ now give
$2\le k\le m-3=6$.  Thus the path representatives $e_k$ and
$e_{k+2}$ used below have valid, distinct indices.

We may now normalize the two endpoint sets without an implicit
preservation assumption. The central relations $k+1\in I_1$ and
$k-1\in I_2$ are exchanged under reversal. In particular,
normalizing the right side does not require a new choice of key
vertices or a relabelling of the missing index after the path is
returned to its original orientation. The full-partition case has
been excluded for every longest path with two high-degree endpoints, and
the endpoint-size argument above gives size at least four for each side.
Thus hypotheses (N1)--(N2) of Lemma~\ref{B:lem:two-sided-normalization}
hold. Apply that lemma.  It either yields
one of the already excluded outcomes or establishes
\[
 I_1=I_1',\qquad I_2=I_2',
\]
with the same single missing index $k$.  In particular,
normalizing the second endpoint does \emph{not} undo primeness
at the first, because its entire four-element prime subset
is preserved by the exchange identity.

It remains to show that $e_1$ and $e_m$ are pair-chargeable. If $e_1$ were not, its third-vertex index cannot lie in $I_2$ by \cref{G:lem:first-edge-crossing}, so it is the unique missing index $k$ and
\[
 e_1=\{v_0,v_1,v_{k+1}\}.
\]
In the left charging, the only possible non-pair contributions are then $e_1$ itself and the path edge whose predecessor is $k$. Since $\deg_H(v_0)\ge7$, at least five of the six $I_1$-pair slots are occupied. \Cref{B:lem:safe-replacement} again forces the full-partition case, contradiction. Hence $e_1$ is pair-chargeable. By reversal, so is $e_m$.

The one-gap charging lemma now applies on each side:
$S=I_1=I_1'$ and $k+1\in I_1$ on the left, and the
reversed hypotheses on the right.  The pair-chargeability of
$e_1,e_m$ was proved in the preceding paragraph.  Since
$\binom{|S|}{2}=6$ and both endpoint degrees are at least seven,
equality holds in the two \emph{injective} charging maps, not
merely in their numerical bounds. Therefore
\begin{equation}\label{B:eq:left-exact}
 \deg_H(v_0)=7,
 \qquad
 e_{k+2}=\{v_0,v_{k+1},v_{k+2}\},
\end{equation}
\begin{equation}\label{B:eq:right-exact}
 \deg_H(v_m)=7,
 \qquad
 e_k=\{v_m,v_{k-1},v_k\},
\end{equation}
and
\[
 \deg_H(v_0,v_{k+1})=\deg_H(v_m,v_k)=1.
\]
Condition~\eqref{B:cond:C} therefore yields
$\deg_H(v_{k+1}),\deg_H(v_k)\ge7$.

Choose an outside edge $b$ witnessing $k-1\in I_2'$ and replace $e_k$ by $b$. This gives the spanning longest Berge path
\begin{equation}\label{B:eq:Q}
 Q=v_ke_{k+1}v_{k+1}e_{k+2}\cdots e_mv_m\,b\,v_{k-1}e_{k-1}\cdots e_1v_0.
\end{equation}
Its representative set is exactly
\[
 (E(P)\setminus\{e_k\})\cup\{b\},
\]
so all representatives are pairwise distinct. Its endpoints $v_k,v_0$ both have degree at least seven. The preceding exclusion of a size-three endpoint applies verbatim to any longest path with two high-degree endpoints; applying it to $Q$ shows that each endpoint set of $Q$ has size at least four. Their endpoint sets are disjoint and $Q$ has nine indices; therefore their total size is eight or nine. If the total is nine, the endpoint sets partition all nine indices and one of them has size four. As in the $q=0$ argument above, a high-degree size-four endpoint is prime, so \cref{G:lem:full-partition} excludes this possibility. Hence both endpoint sets of $Q$ have size four and there is exactly one missing index.

In the displayed order \eqref{B:eq:Q}, write the first vertices as $w_0=v_k$, $w_1=v_{k+1}$, $w_2=v_{k+2}$, and regard $v_0$ as the right endpoint. The second representative is
\[
 e_{k+2}=\{v_0,w_1,w_2\}.
\]
The path uses ten distinct key vertices and the nine distinct representatives $(E(P)\setminus\{e_k\})\cup\{b\}$; this is the full representative ledger, independent of the third vertex of $b$. Notice that the released edge is $e_k$, not $e_{k+2}$;
thus $e_{k+2}$ remains inside $E(Q)$ throughout this exchange.
Thus index $1$ belongs to the right endpoint set of $Q$ by the special
path-edge alternative.  But $\deg_H(v_0,w_1)=1$, and its unique hyperedge
$e_{k+2}$ is already used by $Q$, so index $1$ has no outside witness.

If the right endpoint set were prime, this would already be a contradiction.
Hence it must have the critical $(4,3)$ type, with path-only index $p=1$.
The right-endpoint version of \eqref{B:eq:critical-rel} then forces
$p-1=0$ into the right prime set, contrary to the boundary exclusion
$0\notin I_2(Q)$.  This eliminates the last $n=10$ configuration.
\end{proof}

\subsection{Completion of the boundary-order proof}

\begin{proof}[Proof of \cref{B:thm:main}]
Choose, among all longest Berge paths, one maximizing the sum of its endpoint
degrees.  For $n\in\{9,10\}$, $t=6$.  If both endpoint degrees are at least
$7$, apply \cref{B:thm:high}.  Otherwise one endpoint has degree at most
$6=t$, and Theorem~\ref{G:thm:low} applies with $T=6$. Its degree
hypotheses are exactly \eqref{B:cond:F}--\eqref{B:cond:C}; the shadow is
$2$-connected by assumption. Lemma~\ref{intro:edge-redundancy} shows that the hypotheses imply
that $H$ is $2$-edge-connected. Thus every hypothesis of
Theorem~\ref{G:thm:low} is explicitly accounted for.
\end{proof}

\section{The orders seven and eight}\label{sec:order7}\label{sec:order8}

For $n=7,8$ we reduce the proof to finite enumeration. The reduction
is given below; Appendix~\ref{app:verification} describes the search
and proves completeness. The code uses the Python standard library.

\begin{theorem}[Small orders]\label{S7:thm:main}\label{S8:thm:main}
Let $n\in\{7,8\}$. Every simple $2$-edge-connected $3$-graph $H$ on $n$
vertices with $2$-connected shadow satisfying
\begin{align*}
 d_H(u,v)=2 &\ \Longrightarrow\
 \max\{\deg_H(u),\deg_H(v)\}\ge4,\\
 \deg_H(u,v)=1 &\ \Longrightarrow\
 \min\{\deg_H(u),\deg_H(v)\}\ge4
\end{align*}
has a Hamiltonian Berge cycle.
\end{theorem}

\begin{lemma}[Triangular links at low vertices]\label{small:link}
Suppose that $H$ satisfies the hypotheses of the small-order theorem.
If $\deg_H(v)\le3$, then the non-isolated part of the link of $v$ is a
triangle. Moreover, for $G=\partial_2H$, every vertex has both degrees
at least three, and
\[
 \deg_H(v)=3\quad\Longleftrightarrow\quad \deg_G(v)=3.
\]
\end{lemma}
\begin{proof}
Edge-connectivity gives $\deg_H(v)\ge2$. If $\deg_H(v)\le3$, the
codegree hypothesis makes every non-isolated vertex of its link have
degree at least two. A nonempty simple graph with at most three edges
and minimum positive degree at least two is a triangle. Thus every low
vertex has exactly three shadow neighbours, which form a clique.
Conversely, simplicity gives
$\deg_H(v)\le\binom{\deg_G(v)}2$; if $\deg_G(v)\le3$, this inequality
and the preceding observation force both degrees to equal three.
\end{proof}

The minimum-degree theorem of Kostochka, Luo, and McCourt~\cite{7}
implies Hamiltonicity when $\delta(H)\ge4$: in its Hamiltonian form,
$r=3$ and $\lfloor(n-1)/2\rfloor=3$ for both orders. Consequently a
counterexample has a low vertex. By Lemma~\ref{small:link}, we may label
it $0$ and write
\begin{equation}\label{small:normal}
\begin{gathered}
 V(H)=\{0,1,\ldots,n-1\},\qquad
 E_H(0)=F:=\{012,013,023\},\\
 E(H)=F\cup S,\quad S\subseteq\binom{\{1,\ldots,n-1\}}3.
\end{gathered}
\end{equation}
There are $\binom{n-1}3$ free triples and $\binom{n-1}2-3$ free shadow
pairs. We do not identify isomorphic instances.

\paragraph{Shadow graphs and realizations.}
Enumerate every labelled graph with the six shadow edges of $F$, no
other edge at $0$, and arbitrary remaining pairs. Retain precisely those
passing the following necessary tests:
\begin{enumerate}[label=\textup{(S\arabic*)},leftmargin=*]
\item the graph is $2$-connected;
\item no two vertices of graph degree at most three have distance two;
\item every shadow edge outside $\partial_2F$ lies in a triangle avoiding $0$;
\item the minimum graph degree is at least three, and the neighbourhood
of every degree-three vertex is a clique.
\end{enumerate}
Here (S2) follows from Lemmas~\ref{G:lem:shadow-distance}
and~\ref{small:link}; the other tests follow directly from
\eqref{small:normal} and Lemma~\ref{small:link}.

For a retained graph $G$, let $\mathcal T(G)$ be its triangles avoiding
$0$, viewed as free triples. Let $R(G)$ contain the free triples
$\{v,a,b\}$ for every degree-three vertex $v$ and every pair
$\{a,b\}\subseteq N_G(v)$. The triples in $F$ are already present.

\begin{lemma}[Exact realization conditions]\label{small:realization}
Let $G$ pass \textup{(S1)--(S4)}, and suppose that
$R(G)\subseteq S\subseteq\mathcal T(G)$. Then $F\cup S$ satisfies all
hypotheses of the small-order theorem and has shadow $G$ if and only if:
\begin{enumerate}[label=\textup{(R\arabic*)},leftmargin=*]
\item every edge of $G$ is contained in a triple of $F\cup S$;
\item every vertex of graph degree at least four belongs to at least
four triples of $F\cup S$;
\item every nonempty proper vertex set is crossed by at least two
triples of $F\cup S$.
\end{enumerate}
\end{lemma}
\begin{proof}
Necessity follows from Lemma~\ref{small:link} and the two connectivity
conditions. Conversely, (R1) makes the shadow exactly $G$. At every
degree-three vertex, the forced triples form its entire triangular link,
so its degree is three and each incident shadow pair has codegree two.
All other vertices have hypergraph degree at least four by (R2). Thus
(S2) gives the Fan condition, and a codegree-one pair has two high-degree
endpoints. Finally, (R3) is equivalent to connectivity after deletion of
any one hyperedge, while (S1) gives $2$-connectivity of the shadow.
\end{proof}

\paragraph{Finite verification.}
For each retained shadow, the program searches all realizations satisfying
\textup{(R1)--(R3)} for a Hamilton-free instance. All constraints have
the form $|S\cap M|\ge t$. Starting with $R(G)$, it either rejects a state
already containing a Hamiltonian Berge cycle or branches by adding enough
triples to meet one unmet constraint. Appendix~\ref{app:verification}
proves that every Hamilton-free realization is covered by this recursion.
Cycle detection uses the complete family of inclusion-minimal free-edge
certificates, obtained by enumerating cyclic vertex orders and distinct
representatives. Every retained certificate is also checked by a separate
alternating vertex--hyperedge search.

\begin{table}[htbp]
\centering
\caption{Exhaustive verification in the normal form \eqref{small:normal}.
State counts include initial and rejected states and count a state once
per retained shadow.}\label{small:results}
\begin{tabular}{lrr}
\toprule
Quantity & $n=7$ & $n=8$\\
\midrule
Free triples & 20 & 35\\
Labelled shadows before filtering & 4,096 & 262,144\\
Shadows after (S1)--(S4) & 76 & 10,329\\
Cyclic vertex orders & 72 & 360\\
Minimal cycle certificates & 1,761 & 95,043\\
Search states & 6,824 & 1,778,738\\
Hamilton-free realizations found & 0 & 0\\
\bottomrule
\end{tabular}
\end{table}

\begin{proof}[Proof of the small-order theorem]
A counterexample with minimum degree at least four is excluded by~\cite{7}.
Every other counterexample can be relabelled as in \eqref{small:normal},
and its shadow passes (S1)--(S4). Lemma~\ref{small:realization} places it
in the corresponding realization problem. The complete search described
in Appendix~\ref{app:verification}, with the output in
Table~\ref{small:results}, found no Hamilton-free realization.
\end{proof}

The primary verification is supplied as a Python archive, and the
independent full search as a separate C++20 archive; see
Appendix~\ref{app:verification}.

\begin{proof}[Proof of Theorem~\ref{global:main}]
Lemma~\ref{intro:edge-redundancy} supplies $2$-edge-connectivity.
Theorems~\ref{S7:thm:main}, \ref{B:thm:main}, and
\ref{L:thm:main} then cover $n\in\{7,8\}$, $n\in\{9,10\}$,
and $n\ge11$, respectively.
\end{proof}

\appendix
\section{Specification and completeness of the finite verification}
\label{app:verification}

We describe the finite verification of Theorem~\ref{S7:thm:main}
for $n\in\{7,8\}$, with $F$ as in \eqref{small:normal}. Source code
and reference output are in \path{orders78-verification.zip}.

\subsection{Labelled search space and shadow filters}
Put $W=\{1,\ldots,n-1\}$, $\mathcal T=\binom W3$, and
$m=|\mathcal T|$. Free triples are indexed from $0$ in lexicographic
order. A set $S\subseteq\mathcal T$ is represented by the integer
$s=\sum_{T_i\in S}2^i$. Similarly index the free shadow pairs
$\binom W2\setminus\{12,13,23\}$ and enumerate all their masks in increasing
order. This enumerates every normalized labelled shadow, including
isomorphic copies.

Apply (S1)--(S4) in order: test connectivity after each vertex deletion,
common neighbours of nonadjacent low pairs, triangle support in $W$,
and minimum degree and degree-three neighbourhoods. The successive
survivor counts appear in Table~\ref{app:counts}.

\begin{table}[htbp]
\centering
\caption{Reference counts after successive shadow tests.}\label{app:counts}
\begin{tabular}{crrrr}
\toprule
$n$ & (S1) & (S1)--(S2) & (S1)--(S3) & (S1)--(S4)\\
\midrule
7 & 2,252 & 127 & 124 & 76\\
8 & 168,096 & 17,569 & 15,141 & 10,329\\
\bottomrule
\end{tabular}
\end{table}

For a retained shadow $G$, restrict free variables to $\mathcal T(G)$
and initially select $R(G)$. Triples through $0$ required by a
degree-three neighbourhood are already in $F$. Indeed, if a degree-three
vertex is adjacent to $0$, its other two neighbours must lie in
$N_G(0)=\{1,2,3\}$ by (S4). Thus all forced triples are available.

\subsection{Cycle certificates and a separate witness check}
Fix $0$ first in a cyclic order and write it as
\[
 (0,a,w_1,\ldots,w_{n-3},b),
 \qquad \{a,b\}\in\binom{\{1,2,3\}}2,\quad a<b.
\]
For each choice of $a,b$, enumerate every permutation of
$W\setminus\{a,b\}$. These $3(n-3)!$ orders include every possible
Hamiltonian Berge cycle of \eqref{small:normal}, modulo orientation and
rotation. For every order, enumerate all injective assignments of
representatives from $F\cup\mathcal T$ to its consecutive pairs.
The free representatives used by an assignment form a
\emph{certificate} $C\subseteq\mathcal T$. The two pairs incident with
$0$ use two distinct members of $F$; at most one other member of $F$ can
be used. Hence $|C|\in\{n-3,n-2\}$.

Remove duplicates, then discard a certificate of size $n-2$ precisely
when it contains a certificate of size $n-3$. Let $\mathcal C_n$ be the
remaining family. It is exactly the inclusion-minimal certificate family,
and for every free-edge set $S$,
\begin{equation}\label{app:certificate-equivalence}
 F\cup S\text{ contains a Hamiltonian Berge cycle}
 \quad\Longleftrightarrow\quad
 (\exists C\in\mathcal C_n)\ C\subseteq S.
\end{equation}
For the forward direction, the cycle appears in one enumerated order and
one representative assignment, and its certificate contains a retained
certificate. For the reverse direction, the representative assignment
defining a certificate remains valid in every superset of its free edges.

\begin{table}[htbp]
\centering
\caption{Certificate generation and exhaustive independent witness checks.}
\label{app:certcounts}
\begin{tabular}{lrr}
\toprule
Quantity & $n=7$ & $n=8$\\
\midrule
Distinct certificates of size $n-3$ & 1,266 & 57,348\\
Distinct certificates of size $n-2$, before reduction & 8,115 & 589,923\\
Retained certificates of size $n-2$ & 495 & 37,695\\
Retained certificates checked independently & 1,761 & 95,043\\
Failed witness checks & 0 & 0\\
\bottomrule
\end{tabular}
\end{table}

Store the certificates in a prefix tree of increasing edge indices to
test \eqref{app:certificate-equivalence}. A separate depth-first search
reconstructs $F\cup C$ for every retained mask, alternates vertices and
hyperedges, and checks an explicit length-$n$ cycle with distinct vertices
and representatives. It uses neither the prefix tree nor the generating
assignments. This audits every pruning witness; enumeration completeness
is proved separately here and below.

\subsection{Integer constraints and exhaustive branching}
Write every requirement as $|S\cap M|\ge t$, where
$M\subseteq\mathcal T(G)$. For a pair $uv\in E(G)$ not covered by $F$,
use $M=\{T\in\mathcal T(G):u,v\in T\}$ and $t=1$. For each vertex
of graph degree at least four, use
\[
 M=\{T\in\mathcal T(G):v\in T\},\qquad
 t=4-|\{e\in F:v\in e\}|.
\]
For each nonempty proper vertex set $X$, let $M$ consist of the free
triples meeting both $X$ and its complement, and put
\[
 t=2-|\{e\in F:e\cap X\ne\varnothing,
                    \ e\cap(V\setminus X)\ne\varnothing\}|.
\]
One representative of each complementary pair of cuts suffices. Requirements
with $t\le0$ are vacuous. When several requirements have the same mask,
retain their maximum threshold. Neither simplification changes the feasible sets.

\begin{quote}
\small
\textbf{Search procedure for a retained shadow $G$.}
Begin with $S=R(G)$ and an empty collection of visited states.
\begin{enumerate}[label=\arabic*.,leftmargin=*]
\item If $S$ has been visited for this shadow, return. Otherwise record it.
\item If $S$ contains a member of $\mathcal C_n$, return.
\item For every unmet requirement, calculate its deficit
$d=t-|S\cap M|$ and available set $A=M\setminus S$.
If $|A|<d$, return.
\item If all requirements are met, output $F\cup S$ as a counterexample
and stop with failure status.
\item Choose an unmet requirement minimizing $\binom{|A|}{d}$.
For every $d$-element subset $D$ of $A$, call the procedure with
$S\cup D$.
\end{enumerate}
Tie-breaking and iteration order are deterministic; branches do not
fix edges to be absent.
\end{quote}

\begin{proposition}[Completeness of the search]\label{app:complete}
For each retained shadow, the procedure returns a counterexample if and
only if a Hamilton-free realization satisfying \textup{(R1)--(R3)} exists.
\end{proposition}
\begin{proof}
Any reported set meets the constraints, lies between $R(G)$ and
$\mathcal T(G)$, and contains no cycle certificate. It is therefore a
counterexample by Lemma~\ref{small:realization} and
\eqref{app:certificate-equivalence}.

Conversely, fix such a realization $T$. Initially $R(G)\subseteq T$.
Whenever a visited state satisfies $S\subseteq T$, it cannot contain a
cycle certificate, and each unmet requirement has at least its deficit
$d$ elements available in $T\setminus S$. At least one enumerated child
$S\cup D$ therefore remains a subset of $T$. If all requirements are
already met, $S$ itself is a Hamilton-free realization and is reported.
Otherwise this argument continues while the number of selected triples
strictly increases, so a successful terminal state is reached after at
most $m$ additions along a branch. Reusing a previously visited state
does not remove a possible completion: its outgoing branches depend
only on the state and are searched on its first visit. Thus failure to
find a counterexample excludes every such $T$.
\end{proof}

\subsection{A second exhaustive check at order seven}
The program \path{audit_n7_exhaustive.py} directly tests the original
hypotheses on all $2^{20}=1{,}048{,}576$ normalized hypergraphs. Exactly
661,762 are admissible, and an explicit Hamiltonian Berge cycle was checked
in each. This enumeration shares the normalization and witness checker
with the primary verification, but uses none of its shadow filters,
certificates, prefix tree, or constraint branching.

\subsection{A second exhaustive implementation at both orders}
\label{app:independent-verification}

To check that the finite conclusion does not depend on the Python
certificate generator or its prefix-tree oracle, we also implemented the
normalized verification independently in C++20.  This implementation
constructs the labelled shadow graphs and realization constraints directly,
and uses a depth-first search alternating between key vertices and distinct
representative triples.  Whenever a Hamiltonian Berge cycle is found, its
vertex order and representative indices are reconstructed and independently
checked for pair containment and distinctness before the state is pruned.
No cycle certificates, precomputed witness sets, or prefix trees from the
primary implementation are used.

The second implementation performs the same logically exhaustive deficit
branching, but permits either the first or the last unsatisfied constraint
among those with the smallest number of minimal completions.  Its observed
results are shown below; states are counted once per retained labelled
shadow and need not agree under different tie-breaking rules.
\begin{table}[htbp]
\centering
\caption{Independent C++20 exhaustive search.  No admissible Hamilton-free
realization was found in any run.}\label{app:independent-counts}
\begin{tabular}{lrr}
\toprule
Quantity & $n=7$ & $n=8$\\
\midrule
Shadow graphs enumerated & 4,096 & 262,144\\
Shadows after all four filters & 76 & 10,329\\
States, first minimum-cost requirement & 6,824 & 1,778,738\\
States, last minimum-cost requirement & 4,452 & 1,036,814\\
Admissible Hamilton-free realizations & 0 & 0\\
\bottomrule
\end{tabular}
\end{table}

The recursion is exhaustive for the same elementary reason as
Proposition~\ref{app:complete}: if a feasible Hamilton-free realization
contains the current edge set, some subset of the required cardinality
from any deficient constraint remains inside that realization.  Hence a
child of the recursion preserves containment.  The direct Berge-cycle
search is complete because it enumerates every possible next unused
vertex and every unused representative containing the consecutive pair;
at closure it explicitly checks the final representative from the last
vertex to the first.  It starts at vertex $0$ without loss of generality.
Pruning a state with a verified cycle is safe since cycles persist on
adding hyperedges.

The source code, four outputs, a compiler command, and a reproducibility
check are bundled in
\path{orders78-independent-verification.zip}.  The second program shares
the mathematical normalization and constraint characterization established
above, but uses independent shadow and search code and a different
Hamiltonicity oracle.  In particular, this is an independent full search
of the normalized $n=8$ problem, \emph{not} a direct scan of all
$2^{35}$ free-edge masks.

\subsection{Reproduction and interpretation of the output}
The archive contains the main program \path{verify_orders78.py}, the
separate witness checker \path{audit_certificates.py}, the direct
order-seven audit, and \path{check_reproduction.py} for comparing results.
The README gives file and output specifications; \path{reference_results/}
records the actual execution and environment, and \path{MANIFEST.sha256}
identifies the distributed files. Using Python 3.10 or later, run:
\begin{quote}\small
\begin{verbatim}
python verify_orders78.py --orders 7 8 --outdir results
python audit_n7_exhaustive.py --output results/n7_direct_audit.json
python check_reproduction.py --results results
\end{verbatim}
\end{quote}
These commands regenerate the primary verification with the Python
standard library.  For the independent implementation, a C++20 compiler
(e.g. GCC) is needed; compile and run the program in the independent archive
with the first- and last-tie options as specified in its README.
Certificate digests are SHA-256 hashes of the
increasing distinct masks, each written in unsigned decimal followed by
one ASCII line feed, including the final line, with no header or byte-order
mark. The reference JSON allows exact comparison, apart from run times.

No SAT solver is used. The certificate checker is separate from the
generator, while the completeness argument above covers every normalized
case for $n=7,8$.

\subsection{Computation, independence, and reproducibility limits}
The search exhaustively enumerates labelled instances after the proved
normalization \eqref{small:normal}; it is not a random test, nor does it
assume without proof that every relevant shadow has been sampled.
There are four logically distinct components:
\begin{enumerate}[label=\textup{(V\arabic*)},leftmargin=*]
\item \emph{Reduction:} \cref{small:link,small:realization} establish
that every counterexample of order seven or eight must be represented
by one of the retained shadow/realization instances.
\item \emph{Completeness:} \cref{app:complete} proves that the
constraint-branching recursion cannot discard a Hamilton-free
realization.  The minimal-certificate identity
\eqref{app:certificate-equivalence} justifies the only
Hamiltonicity-based pruning step.
\item \emph{Execution and witness audit:} the archived reference run
reports no counterexample and independently checks an alternating
vertex--edge witness for each retained certificate.  For order seven,
a second program enumerates all $2^{20}$ normalized hypergraphs
without relying on the first program's shadow filters, prefix tree,
or constraint branching.
\item \emph{Independent implementation:} a separate C++20 program
recomputes both finite orders with direct, explicitly checked
vertex--representative cycles, and with two different branching orders.
\end{enumerate}
These checks are independent in the explicitly stated algorithmic
respects, but they are \emph{not} a formal verification of the Python
interpreter or a second implementation of every mathematical reduction.
The direct order-seven audit and the separate C++20 program share the
mathematical normalization, but not the Python certificate oracle.
The C++20 program checks the complete normalized search at order eight;
no direct enumeration of all $2^{35}$ free-edge subsets is claimed.

\paragraph{Reproducibility.}
The two archives contain the Python and C++20 source code, reference
outputs, checked cycle witnesses, and file hashes.  They allow both
searches to be repeated and compared.  They supplement, but do not
replace, the reduction and completeness proof given above.

ChatGPT (OpenAI) assisted in developing the finite-verification code.
The search specification and completeness proof are given above; the
supplementary source and reference outputs permit reproduction of all
reported counts and separate checking of the cycle witnesses.

\subsection*{Acknowledgements}
This work was supported by the National Natural Science Foundation of
China (Grant No.\ 12171088).

\subsection*{Code and data availability}
The supplementary archive \path{orders78-verification.zip} contains the
Python source code, instructions, and reference output for the exhaustive
search at orders $7$ and $8$, including the independent certificate checker
and the direct order-seven verification.  The additional archive
\path{orders78-independent-verification.zip} contains the separate C++20
implementation, four verification outputs, and reproduction instructions.
The primary scripts require only the Python standard library; the second
implementation requires a C++20 compiler. Appendix~\ref{app:verification} describes the
reduction, the search algorithm, and its completeness proof.

\subsection*{Declaration of generative AI and AI-assisted technologies
in the manuscript preparation process}
During the development and preparation of this manuscript, the authors
used ChatGPT (OpenAI) to explore candidate proof arguments, check
representative-hyperedge exchanges, assist with finite-verification code,
and improve the organization and readability of the text. AI-assisted
suggestions were subject to author review and revision. The authors take
full responsibility for the correctness, originality, calculations,
references, and final content of this work. No AI tool is an author.


\begin{thebibliography}{99}
\setlength{\itemsep}{-0.45ex}
\bibitem{1} G. A. Dirac, Some theorems on abstract graphs,
\textit{Proc. London Math. Soc.} (3) 2 (1952), 69--81.

\bibitem{2} O. Ore, Note on Hamilton circuits,
\textit{Amer. Math. Monthly} 67 (1960), 55.

\bibitem{3} L. P\'osa, A theorem concerning Hamilton lines,
\textit{Magyar Tud. Akad. Mat. Kutat\'o Int. K\"ozl.} 7 (1962), 225--226.

\bibitem{4} G. Fan, New sufficient conditions for cycles,
\textit{J. Combin. Theory Ser. B} 37 (1984), 221--227.

\bibitem{5} C. Berge, Nombres de coloration de l'hypergraphe $h$-parti complet,
in \textit{Hypergraph Seminar, Ohio State University 1972},
Lecture Notes in Math. 411, Springer, 1974, 13--20.

\bibitem{6} J.-C. Bermond, A. Germa, M.-C. Heydemann, and D. Sotteau,
Hypergraphes hamiltoniens, in \textit{Probl\`emes combinatoires et th\'eorie
des graphes}, Colloq. Internat. CNRS 260, CNRS, Paris, 1978, 39--43.

\bibitem{7} A. Kostochka, R. Luo, and G. McCourt,
Dirac-type theorems for long Berge cycles in hypergraphs,
\textit{J. Combin. Theory Ser. B} 168 (2024), 159--191.
\href{https://doi.org/10.1016/j.jctb.2024.05.001}{doi:10.1016/j.jctb.2024.05.001}.
For the Hamiltonian formulation used here, see also
\href{https://arxiv.org/abs/2109.12637v3}{arXiv:2109.12637v3}, Theorem~1.7(a).

\bibitem{8} N. Salia, P\'osa-type results for Berge hypergraphs,
\textit{Electron. J. Combin.} 31 (2024), P2.42.
\href{https://doi.org/10.37236/11704}{doi:10.37236/11704}.

\bibitem{9} Y. Li, L. Lu, and R. Luo, An Ore-type theorem for $[3]$-graphs,
\textit{arXiv:2505.12035} (2025).

\bibitem{10} Z. F\"uredi, A. Kostochka, and R. Luo,
Berge cycles in non-uniform hypergraphs,
\textit{Electron. J. Combin.} 27 (2020), Paper No. P3.9.
\href{https://doi.org/10.37236/9346}{doi:10.37236/9346}.

\bibitem{11} L. Lu and Z. Wang,
On Hamiltonian Berge cycles in $[3]$-uniform hypergraphs,
\textit{Discrete Math.} 344 (2021), Paper No. 112462.
\href{https://doi.org/10.1016/j.disc.2021.112462}{doi:10.1016/j.disc.2021.112462}.

\bibitem{12} A. Kostochka, R. Luo, and G. McCourt,
Minimum degree ensuring that a hypergraph is hamiltonian-connected,
\textit{European J. Combin.} 114 (2023), Paper No. 103782.
\href{https://doi.org/10.1016/j.ejc.2023.103782}{doi:10.1016/j.ejc.2023.103782}.

\bibitem{13} J. Zhou and J. Cai,
Ore-type condition for longest cycles in hypergraphs and bipartite graphs,
\textit{Discrete Math.} 349 (2026), Paper No. 115276.
\href{https://doi.org/10.1016/j.disc.2026.115276}{doi:10.1016/j.disc.2026.115276}.

\bibitem{14} M. Coulson and G. Perarnau,
A rainbow Dirac's theorem,
\textit{SIAM J. Discrete Math.} 34 (2020), no. 3, 1670--1692.
\href{https://doi.org/10.1137/18M1218881}{doi:10.1137/18M1218881}.

\bibitem{15} Y. Wang, D. Gerbner, and X. Zhao,
A stability theorem for Berge Hamiltonian cycles under a minimum
degree condition, \textit{arXiv:2608.04603} (2026).
\href{https://arxiv.org/abs/2608.04603}{arXiv:2608.04603}.
\end{thebibliography}
\end{document}